\documentclass[a4paper]{article}
\usepackage{amsthm}
\usepackage{amsmath,amssymb,amscd,amsfonts,amstext,amsxtra}
\usepackage{eufrak,fontenc,fontsmpl,latexsym,%mathpi,
rawfonts}
\theoremstyle{plain}
\newtheorem{theorem}{Theorem}
\newtheorem{lemma}[theorem]{Lemma}
\newtheorem{proposition}[theorem]{Proposition}

\newtheorem{definition}[theorem]{Definition}

\newtheorem*{theorem*}{Theorem}
\newtheorem*{lemma*}{Lemma}
\newtheorem*{proposition*}{Proposition}
\newtheorem*{corollary*}{Corollary}
\newtheorem*{definition*}{Definition}
\newtheorem*{conjecture*}{Conjecture}
\theoremstyle{definition}
\newtheorem*{exercice*}{Exercise}
\newtheorem*{remark*}{Remark}
\newtheorem*{exemple*}{Example}

\newtheorem*{keywords}{Keywords}

\begin{document}
\title{A version of smooth $K$-theory adapted to
the total Chern class}
\author{Alain Berthomieu\\\and Universit\'e de Toulouse, C.U.F.R. J.-F. Champollion\and
I.M.T. (Institut de Math\'ematiques de Toulouse, UMR CNRS n$^\circ$ 5219)\and Campus d'Albi, Place de Verdun,
81012 Albi Cedex, France.
\and{\tt{alain.berthomieu@univ-jfc.fr}}}
\maketitle
\begin{abstract}
A version of smooth $K$-theory is constructed,
which is adapted to the total Chern class instead of the Chern character
(contrarily to the theories considered in \cite{MoiPartie1}, \cite{KaroubiTGCCS}, \cite{KaroubiCCFFHA} and
\cite{BunkeSchick}).

Some total Chern class morphism from this $K$-theory to
Cheeger-Simons differential characters is constructed. This answers
a question raised by U. Bunke \cite{BunkeQuestionOberwolfach}.

Some mistakes (of no consequence) of \cite{MoiSMF} are corrected and some precisions are given about \cite{KaroubiCCFFHA} \S7.19.
\begin{keywords} multiplicative $K$-theory, smooth $K$-theory, Chern-Simons transgression.
\end{keywords}
\noindent{\textbf{AMS-classification:}} Primary: 14F05, 19E20, 57R20, secondary: 14F40, 19D55, 53C05, 55R50.
\end{abstract}
\section{Introduction:}
The name ``Smooth $K$-theory'' refers to an extension of topological $K$-theory
by odd degree differential forms developped by Bunke and Schick
\cite{BunkeSchick}, inspired there by previous work by Hopkins and Singer
\cite{Hopkins}. Related constructions had been previously proposed by Karoubi
\cite{KaroubiTGCCS}, \cite{KaroubiCCFFHA} under the name ``multiplicative
$K$-theory'', which are quotients of subrings of the smooth $K$-theory of \cite{BunkeSchick}.

In all these cases, the Chern character 
plays a crucial role, see
\cite{BunkeSchick} \S1.2.* (particularly the end of \S1.2.4)
for an axiomatic justification of this fact. Consequently, there exists
some Chern character morphisms from these extended $K$-theories to extended
cohomology theories, the prototype of which being the theory of
(rationalised) differential
characters \cite{CheegerSimons} (``multiplicative cohomologies'' in \cite{KaroubiCCFFHA}, which are quotients of subrings of differential characters).

However, the question was raised by U. Bunke in 
\cite{BunkeQuestionOberwolfach} if some counterpart of these
extended $K$-theories corresponding to the total Chern class
(instead of the Chern character) could exist, with some
(nonrationalised) differential character valued total Chern class.
The aim of this paper is to answer positively to Bunke's question: the construction
of the obtained ``multiplicative smooth'' $K$-theory (the word ``multiplicative'' referring to the multiplicativity of the total Chern class)
is presented in definition \ref{flo}, and the total Chern class morphism with values in integral differential characters is obtained from theorem \ref{Angelique} below.
The construction needs to elaborate a theory of multiplicative transgression: this is done in \S\ref{Nadine}, using a modified addition on odd degree differential forms (see \S\ref{modifadd}). The relation between usual and multiplicative transgressions is precised in theorem \ref{maitressePS}.

The drawbacks of the theory exposed here are firstly that only $K^0$-theory
and not higher $K$-theories is constructed, and that the ring structure
(coming from the additivity and the multiplicativity of the Chern character)
is lost (because of the too complicated behaviour of the total Chern class under tensor product of vector bundles). The advantage is to obtain some characteristic class (total Chern class) with values in integral (i.e. nonrationalised) differential characters (in contrast with Chern characters). 

Karoubi gave in \cite{KaroubiCCFFHA} \S7.19 some indications of how to construct a
total Chern class from some of his multiplicative $K$-theories (defined
using the Chern character) with values in some multiplicative cohomology
(proved to be canonically isomorphic to Deligne cohomology).
It is proved here in  theorem \ref{Virginia} that there is some natural map from
usual (``additive'') smooth $K$-theory to multiplicative smooth $K$-theory,
and in theorem \ref{Karen} that
Karoubi's total Chern class construction factorizes through this map.

As a by-product in \cite{MoiSMF}, a first attempt was made in this direction,
with some technical mistakes.
However, there are formal differences between the construction here and
the attempt of \cite{MoiSMF}: they are detailed,
(and the mistakes of \cite{MoiSMF} are corrected) in \S\ref{CorrectionMoi} below.

Finally let's note that a ``multiplicative'' counterpart of the Borel classes of \cite{MoiPartie1} is constucted on multiplicative smooth $K$-theory (in definition \ref{definitiondelaclassedeBorel}) and its relation with the differential character valued total Chern class is explained (see
lemma \ref{decompo} and the last statement of theorem \ref{Angelique}).
\section{Multiplicative transgression:}\label{Nadine}
\subsection{Modified addition on odd degree differential forms:}\label{modifadd}
Let $M$ be some smooth manifold, denote the complex-valued differential forms on $M$ by $\Omega^\bullet(M,{\mathbb C})$, and $\Omega^{\rm{odd}}(M,{\mathbb C})\big/d\Omega^{\rm{even}}(M,{\mathbb C})$ by $\Lambda$.
\begin{definition}
The modified sum of two odd degree differential forms $\alpha$ and $\beta$ is
\[\alpha\stackrel\cap+\beta=\alpha+\beta+\alpha\wedge d\beta\]
\end{definition}
\begin{lemma}
$\stackrel\cap+$ is associative, it induces a commutative law on $\Lambda$. It reduces to ordinary addition if one of the two summands is a cohomology class (a closed form in $\Lambda$).
\end{lemma}
\begin{proof}. The associativity is obvious. The commutativity modulo exact forms too:
\[(\alpha\stackrel\cap+\beta)-(\beta\stackrel\cap+\alpha)=\alpha\wedge d\beta-\beta\wedge d\alpha=d(\alpha\wedge\beta)\]
The last statement is trivial.
\end{proof}
Obviously $\alpha\stackrel\cap+0=0\stackrel\cap+\alpha=\alpha$, and the
opposite element to $\alpha$ is \[\stackrel\cap-\alpha=-\sum_{j=0}^{\infty}(-1)^{j}\alpha\wedge(d\alpha)^{\wedge j}=-\alpha(1+d\alpha)^{-1}\]
(the inverse of $1+d\alpha$ in the sense of power series is a finite sum).

We will need the following formulae (where $\alpha\stackrel\cap-\beta=\alpha\stackrel\cap+(\stackrel\cap-\beta)$)
\begin{equation}\label{obvious}\begin{aligned}
(1+d\alpha)\wedge(1+d\beta)&=1+d(\alpha\stackrel\cap+\beta)\\
(1+d\alpha)\wedge(1+d\beta)^{-1}&=1+d(\alpha\stackrel\cap-\beta)
\end{aligned}
\end{equation}
and the following result : let $F$ and $G$ be the power series defined by
\begin{equation}\label{Christine}F(x)=\frac1x(e^x-1)\qquad{\text{ and }}\qquad G(x)=\frac1x\ln(1+x)\end{equation}
\begin{lemma}\label{labellePatricia} for any elements $\alpha$, $\alpha'$, $\beta$ and $\beta'\in\Lambda$ one has
\begin{align*}(\alpha+\alpha')\wedge F(\alpha+\alpha')&=
\big(\alpha\wedge F(d\alpha)\big)\stackrel\cap+\big(\alpha' \wedge F(d\alpha')\big)\\
\beta\wedge G(d\beta)+\beta'\wedge G(d\beta')&=(\beta\stackrel\cap+\beta')\wedge G\big(d(\beta\stackrel\cap+\beta')\big)
\end{align*}
\end{lemma} 
\begin{proof}
Expand $\theta=\alpha \wedge F(\alpha)\stackrel\cap+\alpha'\wedge F(\alpha')-(\alpha+\alpha')\wedge F(\alpha+\alpha')$, reduce to the same denominator and recombine the obtained terms: this yields
\begin{align*}\theta
&=\frac{d(\alpha\wedge\alpha')}{d\alpha\wedge d\alpha'\wedge(d\alpha+d\alpha')}\wedge\left[d\alpha\wedge\left(e^{d\alpha'}-1\right)-d\alpha'\wedge e^{d\alpha'}\wedge\left(e^{d\alpha}-1\right)\right]\\&=
\frac{d(\alpha\wedge\alpha')}{(d\alpha+d\alpha')}\wedge\left[
F(d\alpha')-e^{d\alpha'}\wedge F(d\alpha)\right]\end{align*}
The fact that $F(x)-e^{x} F(-x)$ vanishes
proves that the two variables power series $F(x)-e^xF(y)$ is divisible by $(x+y)$. Let $H(x,y)$ be the quotient, then $\theta
=d\big(\alpha\wedge\alpha'\wedge H(d\alpha,d\alpha')\big)$ is an exact form. The firt equality is proved (as an equality in $\Lambda$).
The second one is a consequence of the following equivalence
\begin{equation}\label{equivFG}
\beta=\alpha\wedge F(d\alpha)\ \Longleftrightarrow\ \alpha=\beta\wedge G(d\beta)\end{equation}
(it can also be proved in a similar manner as the first one).
\end{proof}
\subsection{Complex conjugation:}\label{complexe}
\begin{definition}
$\theta\in\Omega^{\rm{even}}(M,{\mathbb C})$ is said to be of modulus $1$ if $\theta\wedge\overline\theta=1$.

$\Omega_{[+]}(M,{\mathbb C})$ is the multiplicative (commutative) subgroup of
$\Omega^{\rm{even}}(M,{\mathbb C})$ consisting of even multidegree
differential forms whose degree zero component is equal to $1$ (with
multiplication given by $\wedge$).
\end{definition}
\begin{lemma}\label{SophieMusicale}
Any element in $\Omega_{[+]}(M,{\mathbb C})$ has a unique square root. The square root of a closed form is itself closed. The same holds for real forms.

Any $\sigma\in\Omega_{[+]}(M,{\mathbb C})$ decomposes in a unique way
as a product $\sigma=\rho\wedge\theta$ where $\rho\in\Omega_{[+]}(M,{\mathbb C})$
is a real differential form, and $\theta$ is of modulus $1$.
$\rho$ and $\theta$ are closed if $\sigma$ is closed.

Any element
of $\Omega^{\rm{even}}(M,{\mathbb C})$ of modulus $1$ lies in $\Omega_{[+]}(M,{\mathbb C})$.

A closed even (multi)degree differential form of modulus $1$ whose cohomology class is real, is the sum of the constant $1$ and an exact form.
\end{lemma}
\begin{proof} The square root is obtained from the classical series $\sqrt{1+X}$. Then $\rho=\sqrt{\sigma\wedge\overline\sigma}$ and $\theta=\rho^{-1}\wedge\sigma$. The same results hold for the multiplicative group of even (multi)degree
cohomology whose zero degree part equals $1$. The last assertion is a consequence of the unicity of the decomposition in cohomology.
\end{proof}
The complex conjugation is compatible with $\stackrel\cap+$ and for any $\alpha
\in\Lambda$, $\alpha\stackrel\cap+\overline\alpha$ is real (up to an exact form).
\begin{definition}
If $\gamma\in\Lambda$ verifies $\gamma\stackrel\cap+\overline\gamma=0$ (or equivalently
$\overline\gamma=\stackrel\cap-\gamma$), then $\gamma$ is said to be special imaginary.
\end{definition}
\begin{lemma}\label{deCadix}
For any $\alpha\in\Lambda$, $\alpha\stackrel\cap-\overline\alpha$ is special imaginary.

There exists a special division by $2$ in $\Lambda$
in the sense that for any $\alpha\in\Lambda$, there exists a unique  
$\beta=\genfrac{}{}{0pt}{}{1}{\stackrel\cap2}\alpha\in\Lambda$ such that $\beta\stackrel\cap+\beta=\alpha$.

This special division by $2$ is compatible with complex conjugation.

Any element $\alpha$ of $\Lambda$ uniquely
decomposes as a special sum $\alpha=\beta\stackrel\cap+\gamma$ where $\beta$ is real
and $\gamma$ is special imaginary. (In fact $\beta=
\genfrac{}{}{0pt}{}{1}{\stackrel\cap2}(\alpha\stackrel\cap+\overline\alpha)$ and $\gamma=\genfrac{}{}{0pt}{}{1}{\stackrel\cap2}(\alpha\stackrel\cap-\overline\alpha)$).
\end{lemma}
\noindent With the notations of the lemma, one has $\sqrt{1+d\alpha}=1+d(\genfrac{}{}{0pt}{}{1}{\stackrel\cap2}\alpha)$ and
$\overline\alpha=\beta\stackrel\cap-\gamma$.
\begin{proof}$\beta=\genfrac{}{}{0pt}{}{1}{\stackrel\cap2}\alpha$ is constructed recursively from its degree one part
$\beta^1=\frac12\alpha^1$, then its degree $3$ part $\beta^3=\frac12\left(\alpha^3-\beta^1\wedge d\beta^1\right)$
and so on. The equality of two candidates for being $\genfrac{}{}{0pt}{}{1}{\stackrel\cap2}\alpha$ is proved by reasoning on their lowest degree different component.
The other statements of this lemma are immediate.
\end{proof}
\begin{lemma}\label{sophie}
If $\alpha\in\Lambda$, then $1+d\alpha$ is real if and only if
$\alpha$ is the sum of a closed form and a real form.

If $\alpha\in\Lambda$, then
$1+d\alpha$ is of modulus one if and only if $\alpha$ is the sum
of a closed form and a special imaginary form $\gamma$.
\end{lemma}
\begin{proof}In the first case, the (usual) imaginary part of $\alpha$
is obviously closed. In the second case, decompose $\alpha$ as in the preceding lemma, the conclusion follows from the fact that $\beta\stackrel\cap+\beta$ is closed only if $\beta$ is closed. This is proved inductively on the increasing degree parts of $\beta$.
\end{proof}
Both decompositions (as special sum of a real form and a special imaginary form, or as product of a real form and a modulus one form) exactly correspond to each other by the map
$\alpha\in\Lambda\longmapsto1+d\alpha\in\Omega_{[+]}(M,{\mathbb C})$.
\subsection{Ordinary (``additive'') transgression:}

Let $E$ be a smooth complex vector bundle on $M$. Any connection $\nabla$ on $E$ extends to an
exterior differential operator $d^\nabla$ on smooth
$E$-valued differential forms. The square of $d^\nabla$
is a nondifferential operator, namely the left exterior product
by a degree two ${\rm{End}}E$-valued differential form
known as the ``curvature'' of $\nabla$ and denoted here by
$\nabla^2$.

Chern-Weil theory associates to any conjugation invariant polynomial $P$
on matrices a differential form $P(\nabla^2)$ which is closed, whose
de Rham cohomology class is independent on the choice of the connection
$\nabla$, and which is functorial by pull-back. The polynomials considered here will be the total Chern class $c_{\rm{tot}}(\nabla^2)={\rm{det}}\big({\rm{Id}}+\frac i{2\pi}\nabla^2\big)$ and the Chern character
${\rm{ch}}(\nabla^2)={\rm{Tr}}\big(\exp\big(\frac i{2\pi}\nabla^2\big)\big)$.

Consider $\widetilde M=[0,1]\times M$ and the projection
$p_2\colon[0,1]\times M\longrightarrow M$ on the second factor.
Pull $E$ back on $\widetilde M$ by $p_2$ to obtain a vector
bundle called $\widetilde E$.
For two connections $\nabla_{\!0}$ and $\nabla_{\!1}$ on the
same vector bundle $E$, take any connection $\widetilde\nabla$ on $\widetilde E$ such that its restrictions on $\{0\}\times M$ and
$\{1\}\times M$ verify $\widetilde\nabla\vert_{\{0\}\times M}\cong\nabla_{\!0}$ and $\widetilde\nabla\vert_{\{1\}\times M}\cong\nabla_{\!1}$. The differential form $P(\widetilde\nabla^2)$ is closed, so that its integral along $[0,1]$
\begin{equation}\label{CS}
\widetilde P(\nabla_{\!0},\nabla_{\!1})=\int_{[0,1]}P(\widetilde\nabla^2)
\end{equation}
verifies the following transgression formula
\begin{equation}\label{addtrans}d\widetilde P(\nabla_{\!0},\nabla_{\!1})=P(\nabla^2_{\!1})-P(\nabla^2_{\!0})\end{equation}
Moreover, $\widetilde P(\nabla_{\!0},\nabla_{\!1})$ depends on the choice of $\widetilde\nabla$ only through addition of an exact differential form, it is thus a canonical element of $\Lambda$.

This transgression is additive in the sense that if $\nabla_{\!0}$, $\nabla_{\!1}$ and $\nabla_{\!2}$ are connections on the same vector bundle, then in $\Lambda$ (i.e. up to exact forms):
\begin{equation}\label{Chasles}\widetilde P(\nabla_{\!0},\nabla_{\!2})=\widetilde P(\nabla_{\!0},\nabla_{\!1})+\widetilde P(\nabla_{\!1},\nabla_{\!2})
\end{equation}
In particular, $\widetilde P(\nabla_{\!1},\nabla_{\!0})=-\widetilde P(\nabla_{\!0},\nabla_{\!1})$. And if the polynomial $P$ corresponds to some additive characteristic class (with respect to direct sum of vector bundles, as the Chern character for instance), then
\begin{equation}\label{directsum}\widetilde P(\nabla_{\!E,0}\oplus\nabla_{\!F,0},\nabla_{\!E,1}\oplus\nabla_{\!F,1})=\widetilde P(\nabla_{\!E,0},\nabla_{\!E,1})+\widetilde P(\nabla_{\!F,0},\nabla_{\!F,1})
\end{equation}
if $(\nabla_{\!E,0},\nabla_{\!E,1})$ and $(\nabla_{\!F,0},\nabla_{\!F,1})$ are two couples of connections on the vector bundles $E$ and $F$ on $M$ respectively, and $\nabla_{\!E,0}\oplus\nabla_{\!F,0}$ and $\nabla_{\!E,1}\oplus\nabla_{\!F,1}$ denote the associated direct sum connections on $E\oplus F$. Moreover
\begin{lemma}\label{wedge}
If $P$ and $Q$ are two conjugation invariant polynomials
on matrices, then $P\wedge Q$ is one too, and one has the following equality in $\Lambda$:
\begin{align*}\widetilde{P\wedge Q}(\nabla_{\!0},\nabla_{\!1})&=
\widetilde P(\nabla_{\!0},\nabla_{\!1})\wedge Q(\nabla_{\!0}^2)
+P(\nabla_{\!1}^2)\wedge\widetilde Q(\nabla_{\!0}, \nabla_{\!1})
\\&=P(\nabla_{\!0}^2)\wedge\widetilde Q(\nabla_{\!0},\nabla_{\!1})
+\widetilde P(\nabla_{\!0},\nabla_{\!1})\wedge Q(\nabla_{\!1}^2)
\end{align*}
\end{lemma}
\begin{proof}
The two right hand sides differ from $d\big(\widetilde P(\nabla_{\!0},\nabla_{\!1})\wedge\widetilde Q(\nabla_{\!0}, \nabla_{\!1})\big)$
so that they define equal elements of $\Lambda$.

Consider $E\oplus E$, endowed with a direct sum connection, and the
polynomial $P\wedge Q$, where $P$ is evaluated on the curvature of the connection of the first copy of $E$ in the direct sum $E\oplus E$, and $Q$ on the second one.

Of course $P\wedge Q$ is not conjugation invariant on ${\rm{End}}
(E\oplus E)$, but it is invariant under conjugation by elements
of ${\rm{End}}E\oplus{\rm{End}}E$, i.e. endomorphisms of $E\oplus E$ which respect the decomposition. This means that it can be evaluated on the curvature of any
direct sum connection (which does not need to be the direct sum of twice the same connection on $E$).

With this restriction, the preceding theory works and there is a notion of
$\widetilde{P\wedge Q}$ in this context. The equality above then follows from the possibility to calculate this $\widetilde{P\wedge Q}(\nabla_{\!0}\oplus\nabla_{\!0},
\nabla_{\!1}\oplus\nabla_{\!1})$ with different connections on the bundle $\widetilde{E\oplus
E}=\widetilde E\oplus\widetilde E$ over $[0,1]\times M$. If one first changes the connection on the first copy of $E$ and then on the second copy, one obtains the
first right hand side, if one puts on $\widetilde{E\oplus E}=\widetilde E\oplus\widetilde E$ the direct sum of
twice the same connection on each copy of $\widetilde E$, one obtains a differential form which is easily seen to equal $\widetilde{P\wedge Q}(\nabla_{\!0},\nabla_{\!1})$ in the sense of the first definition of $\widetilde{P\wedge Q}$.
\end{proof}
Finally, Chern-Simons transgression forms are locally gauge invariant:
if $g\colon E\overset\sim\longrightarrow E$ is a global smooth bundle automorphsm which is isotopic to the identity of $E$ (i.e. smoothly homotopic throw isomorphisms), then
$\widetilde P(\nabla,g^*\nabla)$ is an exact form; precisely, 
if $g_t$ is a smooth family of automorphisms of $E$
\begin{equation}\label{[}\frac d{dt}\widetilde P(\nabla,g_t^*\nabla)=d\left(\frac\partial
{\partial b}P\left((g_t^*\nabla)^2+bg_t^{-1}\frac{\partial g_t}{\partial t}\right)\Big\vert_{b=0}\right)
\end{equation}
\subsection{``Multiplicative'' transgression:}
This kind of transgression is adapted to conjugation invariant polynomials which are multiplicative with respect to direct sums and have a zero degree component equal to $1$ (e.g. $c_{\rm{tot}}$). Such a polynomial $P$ can be inverted (in power series).
\begin{definition}
For two connections $\nabla_{\!0}$ and $\nabla_{\!1}$ on the same
vector bundle $E$, put
\[\widehat P(\nabla_{\!0},\nabla_{\!1})=
\widetilde P(\nabla_{\!0},\nabla_{\!1})\wedge P^{-1}(\nabla_{\!0}^2)\]
\end{definition}
$P^{-1}(\nabla_{\!0}^2)$ being closed, $\widehat P(\nabla_{\!0},\nabla_{\!1})$ is well defined in $\Lambda$. 
\begin{lemma} $\widehat P$ verifies the following ``multiplicative'' transgression formula:
\begin{equation}\label{multrans}
1+d\widehat P(\nabla_{\!0},\nabla_{\!1})=P(\nabla_{\!1}^2)\wedge P(\nabla_{\!0}^2)^{-1}\end{equation}
If $\nabla_0$, $\nabla_1$ and $\nabla_2$ are three connections on $E$, then
\begin{equation}\label{multcocycle}
\widehat P(\nabla_{\!0},\nabla_{\!2})=\widehat P(\nabla_{\!1},\nabla_{\!2})\stackrel\cap+\widehat P(\nabla_{\!0},\nabla_{\!1})
\end{equation}
\end{lemma}
\begin{proof}\eqref{multrans} is a direct consequence of \eqref{addtrans}. \eqref{multcocycle} is a consequence of \eqref{Chasles}.\end{proof}
The following consequence of \eqref{multcocycle}:
\begin{equation}\label{echange}
\widehat P(\nabla_{\!1},\nabla_{\!0})=\ \stackrel\cap-
\widehat P(\nabla_{\!0},\nabla_{\!1})
\end{equation}
should replace the second equality in formula (9) of \cite{MoiSMF} (which is false there).
\begin{lemma} One has the following formula:
\[\widehat{P^{-1}}(\nabla_{\!0},\nabla_{\!1})=\widehat P(\nabla_{\!1},
\nabla_{\!0})\]
\end{lemma}
It was already stated as first equality in formula (9) of \cite{MoiSMF} (where it is true).
\begin{proof}
If $P_{[>0]}$ stands for the positive degree part of $P$ then $P=1+P_{[>0]}$ and
\[P^{-1}=1+\sum_{i=1}^\infty(-1)^i\big(P_{[>0]}\big)^{\wedge i}\]
so that following lemma \ref{wedge}
\begin{align*}
\widetilde{P^{-1}}(\nabla_{\!0},\nabla_{\!1})&=
\sum_{i=1}^\infty(-1)^i\sum_{j=0}^i\big(P_{[>0]}(\nabla_{\!0}^2)\big)^{\wedge j}
\wedge\widetilde P(\nabla_{\!0},\nabla_{\!1})\wedge
\big(P_{[>0]}(\nabla_{\!1}^2)\big)^{\wedge(i-j+1)}\\
=-\Big(1+\sum_{j=1}^\infty&(-1)^j\big(P_{[>0]}(\nabla_{\!0}^2)\big)^{\wedge j}\Big)
\wedge\widetilde P(\nabla_{\!0},\nabla_{\!1})\wedge
\Big(1+\sum_{k=1}^\infty(-1)^k\big(P_{[>0]}(\nabla_{\!1}^2)\big)^{\wedge k}\Big)\\
=-P^{-1}(\nabla_{\!0}^2&)\wedge\widetilde P(\nabla_{\!0},\nabla_{\!1})\wedge
P^{-1}(\nabla_{\!1}^2)
\end{align*}
The desired relation follows.
\end{proof}

For a direct sum $E\oplus F$, if $P$ is multiplicative (with respect to direct sums as the total Chern class) then so is $P^{-1}$ and (with notations as in \eqref{directsum}):
\begin{equation}\label{multsum}
\begin{aligned}\widehat P(\nabla_{\!E,0}&\oplus\nabla_{\!F,0},\nabla_{\!E,1}\oplus\nabla_{\!F,1})=\\&=
\widetilde P(\nabla_{\!E,0},\nabla_{\!E,1})\wedge P\big((\nabla_{\!F,0})^2\big)\wedge P^{-1}\big((\nabla_{\!E,0}\oplus\nabla_{\!F,0})^2\big)+\\&\qquad+P\big((\nabla_{\!E,1})^2\big)\wedge\widetilde P(\nabla_{\!F,0},\nabla_{\!F,1})\wedge P^{-1}\big((\nabla_{\!E,0}\oplus\nabla_{\!F,0})
^2\big)\\
&=\widehat P(\nabla_{\!E,0},\nabla_{\!E,1})+\widehat P(\nabla_{\!F,0},\nabla_{\!F,1})\wedge\big(1+d\widehat P(\nabla_{\!E,0},\nabla_{\!E,1})\big)\\
&=\widehat P(\nabla_{\!F,0},\nabla_{\!F,1})\stackrel\cap+
\widehat P(\nabla_{\!E,0},\nabla_{\!E,1})
\end{aligned}
\end{equation}
The local gauge invariance \eqref{[} also holds for multiplicative transgression forms, namely if $g\colon E\overset\sim\longrightarrow E$ is a global smooth bundle automorphsm which is isotopic to the identity of $E$, then
$\widehat P(\nabla,g^*\nabla)$ is an exact form.

Finally, one has the following relation between additive and multiplicative transgression: let $\phi\colon
\Omega(M,{\mathbb C})\longrightarrow\Omega(M,{\mathbb C})$
be the additive group morphism which multiplies degree $2k$ and
degree $2k-1$ differential forms by $(-1)^{k-1}(k-1)!$, and degree
$0$ forms by $0$. Then
\begin{theorem}\label{maitressePS} For any connections $\nabla_{\!0}$ and $\nabla_{\!1}$ on some same vector bundle
\[\widehat c_{\rm{tot}}(\nabla_{\!0},\nabla_{\!1})=
\phi\big(\widetilde{\rm{ch}}(\nabla_{\!0},\nabla_{\!1})\big)\wedge
F\big[\phi\big(d\widetilde{\rm{ch}}(\nabla_{\!0},\nabla_{\!1})\big)\big]\]
\end{theorem}
Using the fact that $d\circ\phi=\phi\circ d$ when restricted to odd degree forms, the inverse formula follows \eqref{equivFG}:
\[\widetilde{\rm{ch}}(\nabla_{\!0},\nabla_{\!1})
=\phi^{-1}\big[\widehat c_{\rm{tot}}(\nabla_{\!0},\nabla_{\!1})
\wedge G\big(d\widehat c_{\rm{tot}}(\nabla_{\!0},\nabla_{\!1})\big)\big]\]
\begin{proof}
Consider some connection $\widetilde\nabla$ on the pull-back of the bundle on $M\times[0,1]$ with expected restrictions $\nabla_{\!0}$
and $\nabla_{\!1}$ on $M\times\{0\}$ and $M\times\{1\}$. Call
$\nabla_{\!t}$ the restriction to $M\times\{t\}$. We will
prove the equality for the couple $(\nabla_{\!0},\nabla_{\!t})$
for all $t$ by deriving it respect to $t$. Of course the equality trivialy hols for $t=0$.

Let ${\rm{ch}}(\widetilde\nabla^2)=\alpha_t+dt\wedge\beta_t$
where $\alpha_t={\rm{ch}}(\nabla_{\!t}^2)$ and $\beta_t$ are $t$-depending elements of $\Omega(M,{\mathbb C})$. In one hand, $\frac d{dt}\widetilde{\rm{ch}}(\nabla_{\!0},\nabla_{\!t})=\beta_t$ so that
\begin{align*}
\frac d{dt}&\left[\phi\big(\widetilde{\rm{ch}}(\nabla_{\!0},\nabla_{\!t})\big)\wedge
F\big[\phi\big(d\widetilde{\rm{ch}}(\nabla_{\!0},\nabla_{\!t})\big)\big]\right]=\\&=\phi(\beta_t)\wedge
F\big[\phi\big(d\widetilde{\rm{ch}}(\nabla_{\!0},\nabla_{\!t})\big)\big]+\phi\big(\widetilde{\rm{ch}}(\nabla_{\!0},\nabla_{\!t})\big)\wedge
F'\big[\phi\big(d\widetilde{\rm{ch}}(\nabla_{\!0},\nabla_{\!t})\big)\big]\wedge \phi(d\beta_t)
\end{align*}
where $F'$ is the derived power series of $F$. Up to an exact form, this equals
\begin{align*}\ &=\phi(\beta_t)\wedge\left[
F\big[\phi\big(d\widetilde{\rm{ch}}(\nabla_{\!0},\nabla_{\!t})\big)\big]+\phi\big(d\widetilde{\rm{ch}}(\nabla_{\!0},\nabla_{\!t})\big)\wedge
F'\big[\phi\big(d\widetilde{\rm{ch}}(\nabla_{\!0},\nabla_{\!t})\big)\big]\right]\\
&=\phi(\beta_t)\wedge\exp\left(\phi\big(d\widetilde{\rm{ch}}(\nabla_{\!0},\nabla_{\!t})\big)\right)
\end{align*}
since $F(x)+xF'(x)=e^x$. In the other hand, it was proved in
\cite{MoiSMF} lemma 4.1 that $c_{\rm{tot}}(\nabla^2)=\exp\big(\phi({\rm{ch}}(\nabla^2))\big)$ (for any connection $\nabla$ on any bundle) so that
\begin{align*}c_{\rm{tot}}(\widetilde\nabla^2)&=\exp\big(\phi(a_t+dt\wedge\beta_t)\big)=\exp\big(\phi(\alpha_t)\big)\wedge
\big(1+dt\wedge\phi(\beta_t)\big)\\
\frac d{dt}\widehat c_{\rm{tot}}(\nabla_{\!0},\nabla_{\!t})&=
c_{\rm{tot}}^{-1}(\nabla_{\!0}^2)\wedge\exp\big(\phi(\alpha_t)\big)\wedge\phi(\beta_t)=\exp\big(\phi(\alpha_t)-\phi(\alpha_0)\big)\wedge\phi(\beta_t)
\end{align*}
The desired equality follows the obvious relation $\alpha_t-\alpha_0=d\widetilde{\rm{ch}}(\nabla_{\!0},\nabla_{\!t})$.
\end{proof}
\section{Smooth $K$-theory adapted to $c_{\rm{tot}}$:}
\subsection{``Multiplicative'' smooth $K^0$-theory:}\label{juriste!}
Consider some triple $(E,\nabla_{\!E},\alpha)$ where $\nabla_{\!E}$ is a connection on the complex vector bundle $E$ over $M$, and $\alpha\in\Lambda$. If there is some
vector bundle isomorphism $f\colon E\overset\sim\longrightarrow F$ and if $\nabla_{\!F}$ is a connection on $F$, then the two following triples will be considered as equivalent:
\begin{equation}\label{reldeq}
(E,\nabla_{\!E},\alpha)=\big(F,\nabla_{\!F},\widehat c_{\rm{tot}}(\nabla_{\!E},f^*\nabla_{\!F})\stackrel\cap+\alpha\big)
\end{equation}

In view of \eqref{multcocycle}, this is intended to ensure the transitivity of this relation, namely if $g\colon F\overset\sim\longrightarrow G$ is another vector bundle isomorphism (and $\nabla_{\!G}$ is any connection on the vector bundle $G$ over $M$), then
\begin{align*}\big(G,\nabla_{\!G},\widehat c_{\rm{tot}}\big(\nabla_{\!E},&(g\circ f)^*\nabla_{\!G}\big)\stackrel\cap+\alpha\big)=\\&\qquad=\big(G,\nabla_{\!G},\widehat c_{\rm{tot}}(\nabla_{\!F},g^*\nabla_{\!G})\overset\cap+\widehat c_{\rm{tot}}(\nabla_{\!E},
f^*\nabla_{\!F})\overset\cap+\alpha\big)
\end{align*}
In the same way, the symmetry is ensured by \eqref{echange}:
\begin{align*}\big(F,\nabla_{\!F},\widehat c_{\rm{tot}}&(\nabla_{\!E},f^*\nabla_{\!F})\stackrel\cap+\alpha\big)
=\\&=\big(E,\nabla_{\!E},\widehat c_{\rm{tot}}(f^*\nabla_{\!F},\nabla_{\!E})\stackrel\cap+\widehat c_{\rm{tot}}(\nabla_{\!E},f^*\nabla_{\!F})\stackrel\cap+\alpha\big)=(E,\nabla_{\!E},\alpha)
\end{align*}
\begin{definition}\label{flo}
The ``multiplicative'' smooth $K$-theory group $\overset*K{}^0_{c_{\rm{tot}}}(M)$ is the quotient of the free abelian group generated by equivalence classes of triples as above by the following relation (for any vector bundle $H$ on $M$ with connection $\nabla_{\!H}$):
\begin{equation}\label{relsum}
(E,\nabla_{\!E},\alpha)+(H,\nabla_{\!H},\beta)=\big(E\oplus H,\nabla_{\!E}\oplus\nabla_{\!H},
\alpha\stackrel\cap+\beta\big)
\end{equation}
\end{definition}

If as above $f\colon E\overset\sim\longrightarrow F$ and $h\colon H\overset\sim\longrightarrow K$
are vector bundle isomorphisms, then for any connections $\nabla_{\!H}$ and $\nabla_{\!K}$ on $H$ and $K$:
\begin{align*}
\Big(F&\oplus K,\nabla_{\!F}\oplus\nabla_{\!K},\widehat c_{\rm{tot}}\big(\nabla_{\!E}\oplus\nabla_{\!H},
(f\oplus h)^*(\nabla_{\!F}\oplus\nabla_{\!K})\big)\stackrel\cap+(\alpha\stackrel\cap+\beta)\Big)=\\
&=\Big(F\oplus K,\nabla_{\!F}\oplus\nabla_{\!K},\big(\widehat c_{\rm{tot}}(\nabla_{\!E},
f^*\nabla_{\!F})\stackrel\cap+\alpha\big)\stackrel\cap+\big(\widehat c_{\rm{tot}}(\nabla_{\!H},
h^*\nabla_{\!K})\stackrel\cap+\beta\big)\Big)\\
&=\big(F,\nabla_{\!F},\widehat c_{\rm{tot}}(\nabla_{\!E},
f^*\nabla_{\!F})\stackrel\cap+\alpha\big)+\big(K,\nabla_{\!K},\widehat c_{\rm{tot}}(\nabla_{\!H},
h^*\nabla_{\!K})\stackrel\cap+\beta\big)
\end{align*}
(because of \eqref{multsum}). This proves the compatibility of relations \eqref{reldeq} and \eqref{relsum}.

The functoriality of $\widehat c_{\rm{tot}}$ by pullbacks makes $\overset*K{}^0_{c_{\rm{tot}}}$ to a contravariant functor from smooth manifolds to commutative groups. 
\subsection{Exact sequence:}
\begin{definition} For $\gamma\in\Lambda$ and $(E,\nabla_{\!E},\alpha)$ one defines
\[\iota(\gamma)=(E,\nabla_{\!E},\gamma\stackrel\cap+\alpha)
-(E,\nabla_{\!E},\alpha)\in\overset*K{}^0_{c_{\rm{tot}}}(M)\]
\end{definition}
\begin{lemma}
$\iota(\gamma)$ does not depend on the triple $(E,\nabla_{\!E},\alpha)$ used to
compute it. $\gamma\longmapsto \iota(\gamma)$ is a morphism from $\Lambda$ endowed with the group law $\stackrel\cap+$
to $\overset*K{}^0_{c_{\rm{tot}}}(M)$.
\end{lemma}
\begin{proof}
Take any two elements $(E,\nabla_{\!E},\alpha)$
and $(H,\nabla_{\!H},\beta)$ of $\overset*K{}^0_{c_{\rm{tot}}}(M)$, then
\begin{align*}
\big((E,&\nabla_{\!E},\gamma\stackrel\cap+\alpha)-(E,\nabla_{\!E},\alpha)\big)
-\big((H,\nabla_{\!H},\gamma\stackrel\cap+\beta)-(H,\nabla_{\!H},\beta)\big)=\\
&=\big((E,\nabla_{\!E},\gamma\stackrel\cap+\alpha)+(H,\nabla_{\!H},\beta)\big)
-\big((H,\nabla_{\!H},\gamma\stackrel\cap+\beta)+(E,\nabla_{\!E},\alpha)\big)\\
&=(E\oplus H,\nabla_{\!E}\oplus\nabla_{\!H},\gamma\stackrel\cap+\alpha\stackrel\cap+\beta)
-(H\oplus E,\nabla_{\!H}\oplus\nabla_{\!E},\gamma\stackrel\cap+\beta\stackrel\cap+\alpha)=0
\end{align*}
The second statement of the lemma is obvious.
\end{proof}

An element of $K^1_{\rm{top}}(M)$ can be represented by some vector bundle
$\xi$ on $S^1\times M$ whose restrictions to $\{pt\}\times M$ are topologically trivial
for any point $pt$ of $S^1$. Choose any connection on $\xi$ whose restriction to
$\{1\}\times M$ is the canonical connection $d$ on the trivial vector bundle
$\xi\vert_{\{1\}\times M}$. The parallel transport along $S^1$ provides a global vector bundle automorphism $\psi$ of $\xi\vert_{\{1\}\times M}$. The (suspended) total Chern
class of $\xi$ is the odd degree cohomology class defined by $\widetilde c_{\rm{tot}}
(d,\psi^*d)$ or $\widehat c_{\rm{tot}}
(d,\psi^*d)$ indifferently. This provides a morphism $K^1_{\rm{top}}(M)\longmapsto
H^{\rm{odd}}(M,{\mathbb C})$ which can be extended to take its values in $\Lambda$,
(the restriction to $H^{\rm{odd}}(M,{\mathbb C})$ of the two group laws $+$ and $\stackrel\cap+$ coincide).

Note however that if one adds to $\xi$ a vector bundle $\eta$ obtained by pullback from $M$ (by the projection of $S^1\times M$ onto the second factor)
with pullback connection $\nabla_{\!\eta}$, then from lemma \ref{wedge}, one obtains in $\Lambda$:
\begin{equation}\label{Patricia}\widetilde c_{\rm{tot}}\big(d\oplus\nabla_{\!\eta},(\psi\oplus{\rm{Id}}_\eta)^*(d\oplus\nabla_{\!\eta})\big)=\widetilde c_{\rm{tot}}(d,\psi^*d)\wedge c_{\rm{tot}}(\nabla_{\!\eta}^2)\end{equation}
while from \eqref{multcocycle} and \eqref{multsum}, one has (in $\Lambda$) for any connection $\nabla$ on $\xi\vert_{\{1\}\times M}$
\begin{equation}\label{K1Sctot}
\begin{aligned}
\widehat c_{\rm{tot}}&\big(d\oplus\nabla_{\!\eta},(\psi\oplus{\rm{Id}}_\eta)^*(d\oplus\nabla_{\!\eta})\big)=\widehat c_{\rm{tot}}(d,\psi^*d)\stackrel\cap+ \widehat c_{\rm{tot}}(\nabla_{\!\eta},\nabla_{\!\eta})=\widehat c_{\rm{tot}}(d,\psi^*d)=\\&=\widehat c_{\rm{tot}}(d,\nabla)\stackrel\cap+\widehat c_{\rm{tot}}(\nabla,\psi^*\nabla)\stackrel\cap+\widehat c_{\rm{tot}}(\psi^*\nabla,\psi^*d)=\widehat c_{\rm{tot}}(\nabla,\psi^*\nabla)
\end{aligned}\end{equation}
Thus $\widehat c_{\rm{tot}}$
is much better adapted than $\widetilde c_{\rm{tot}}$ to represent the suspended total Chern class on $K^1_{\rm{top}}(M)$ (as was already remarked in formula (15) of \cite{MoiSMF}).

Consider the forgetful map  ${\mathcal F}\colon(E,\nabla,\alpha)\in\overset*K{}^0_{c_{\rm{tot}}}(M)\longmapsto[E]\in K^0_{\rm{top}}(M)$ (with value in the topological $K$-theory of $M$).

\begin{proposition}
The following sequence is exact
\begin{equation}\label{pasinvarhomot}
K^1_{\rm{top}}(M)\overset{Sc_{\rm{tot}}}\longrightarrow \Lambda\overset\iota
\longrightarrow\overset*K{}^0_{c_{\rm{tot}}}(M)\overset{\mathcal F}
\longrightarrow K^0_{\rm{top}}(M)\longrightarrow\{0\}
\end{equation}
\end{proposition}
\begin{proof}
The surjectivity of ${\mathcal F}$ is tautological.

The nullity of ${\mathcal F}\circ\iota$ is tautological.

From \eqref{K1Sctot}, if $\alpha\in\Lambda$ lies in the image of $Sc_{\rm{tot}}$, then $\alpha=\widehat c_{\rm{tot}}(\nabla,\psi^*\nabla)$ for some connection $\nabla$ on the trivial vector bundle on $M$ as above
(there the trivial vector bundle was $\xi\vert_{\{1\}\times M}$).
Denote this trivial bundle by ${\mathbb C}^{{\rm{rk}}\xi}$, then
\[\iota(\alpha)=\big({\mathbb C}^{{\rm{rk}}\xi},\nabla,\widehat c_{\rm{tot}}(\nabla,\psi^*\nabla)\big)-({\mathbb C}^{{\rm{rk}}\xi},\nabla,0)\]
vanishes because of \eqref{reldeq}. Thus $\iota\circ Sc_{\rm{tot}}$ vanishes.

Any element of $\overset*K{}^0_{c_{\rm{tot}}}(M)$ can be represented by a difference
of the form $(E,\nabla_{\!E},\alpha)-(H,\nabla_{\!H},\beta)$. If $[E]-[H]=0\in K^0_{\rm{top}}(M)$, then there exists some vector bundle $L$ such that $E\oplus L$ and
$H\oplus L$ are isomorphic.
Let $\nabla_{\!L}$ be any connection on $L$ and
let $\varphi\colon E\oplus L\overset\sim\longrightarrow
H\oplus L$ be such an isomorphism. Then one has the following equalities in $\overset*K{}^0_{c_{\rm{tot}}}(M)$:
\begin{align*}(E,\nabla_{\!E},\alpha)&-(H,\nabla_{\!H},\beta)=(E\oplus L,\nabla_{\!E}\oplus\nabla_{\!L},\alpha)-(H\oplus L,\nabla_{\!H}\oplus\nabla_{\!L},\beta)\\
&=\Big(H\oplus L,\nabla_{\!H}\oplus\nabla_{\!L},\widehat c_{\rm{tot}}\big(\nabla_{\!E}\oplus\nabla_{\!L},\varphi^*(\nabla_{\!H}\oplus\nabla_{\!L})\big)\stackrel\cap+\alpha\Big)\\&\qquad\qquad\qquad\qquad\qquad\qquad\qquad\qquad\, -(H\oplus L,\nabla_{\!H}\oplus\nabla_{\!L},\beta)\\&=\iota\, \Big(\widehat c_{\rm{tot}}\big(\nabla_{\!E}\oplus\nabla_{\!L},\varphi^*(\nabla_{\!H}\oplus\nabla_{\!L})\big)\stackrel\cap+\alpha\stackrel\cap-\beta\Big)
\end{align*}
which proves the exactness at $\overset*K{}^0_{c_{\rm{tot}}}(M)$.

If $(E,\nabla_{\!E},\alpha)=(F,\nabla_{\!F},\beta)\in\overset*K{}^0_{c_{\rm{tot}}}(M)$, then there exists some vector bundle $G$ on $M$, endowed with some connection $\nabla_{\!G}$, and a vector bundle isomorphism $\varphi\colon E\oplus G\overset\sim\longrightarrow F\oplus G$ such that $\beta=\alpha\stackrel\cap+\widehat c_{\rm{tot}}\big(\nabla_{\!E}\oplus\nabla_{\!G},\varphi^*(\nabla_{\!F}\oplus\nabla_{\!G})\big)$.

Thus if $\iota(\alpha)=(E,\nabla_{\!E},\alpha)-(E,\nabla_{\!E},0)$ vanishes in $\overset*K{}^0_{c_{\rm{tot}}}(M)$, then there exists $G$ on $M$, with connection $\nabla_{\!G}$, and an isomorphism $\varphi\colon E\oplus G\overset\sim\longrightarrow E\oplus G$ such that $\alpha=\widehat c_{\rm{tot}}\big(\nabla_{\!E}\oplus\nabla_{\!G},\varphi^*(\nabla_{\!E}\oplus\nabla_{\!G})\big)$.

Choose some vector bundle $H$ such that $E\oplus G\oplus H$ is trivial, then clearly from \eqref{K1Sctot}, $\alpha$ is the image by $Sc_{\rm{tot}}$ of the
element of $K^1_{\rm{top}}(M)$ represented by the global bundle automorphism $\varphi\oplus{\rm{Id}}_H$ of the trivial vector bundle $E\oplus G\oplus H$.
\end{proof}
\subsection{Differential form valued total Chern class on $\overset*K{}^0_{c_{\rm{tot}}}(M)$:}\label{jorba}
\begin{definition}\label{italiénisante}
For any element $(E,\nabla,\alpha)$ of $\overset*K{}^0_{c_{\rm{tot}}}(M)$, put
\[\overset*c_{\rm{tot}}(E,\nabla,\alpha)=c_{\rm{tot}}(\nabla^2)\wedge(1+d\alpha)^{-1}\]
\end{definition}
It follows from \eqref{obvious}, \eqref{reldeq} and \eqref{multrans}
that this defines a map from $\overset*K{}^0_{c_{\rm{tot}}}(M)$ to
$\Omega^{\rm{even}}(M,{\mathbb C})$ (with no need of quotienting by exact forms).
Moreover, it follows from \eqref{obvious} and \eqref{relsum} that this map is a morphism if $\Omega(M,{\mathbb C})$ is endowed
with its multiplicative structure given by $\wedge$.

Note also the following tautological properties:
\begin{align}
c_{\rm{tot}}\big({\mathcal F}(E,\nabla,\alpha)\big)&=c_{\rm{tot}}([E])=\big[\overset*c_{\rm{tot}}(E,
\nabla,\alpha)\big]\, \in\, H^{\rm{even}}(M,{\mathbb C})\\
\overset*c_{\rm{tot}}\big(-\iota(\gamma)\big)&=1+d\gamma
\end{align}
They could be used in some attempt to axiomatize ``multiplicative'' smooth $K$-theory in the same way
as the corresponding ``additive'' theory is in \cite{BunkeSchick} (see \S\S1.2.1 and 2.4.3 in \cite{BunkeSchick}).

The subgroup ${\mathcal M}K^0_{c_{\rm{tot}}}(M)$ of $\overset*K{}^0_{c_{\rm{tot}}}(M)$ consisting of elements whose total Chern class is equal to $1$ is something like a ``nonfree
multiplicative'' $K$-theory group: it is the multiplicative transgression counterpart of Karoubi's multiplicative $K$-theory group corresponding to flat vector bundles \cite{KaroubiTGCCS} EXEMPLE 3 of the introduction and \cite{KaroubiCCFFHA} \S5.5.
(A group of the same kind was considered in \cite{MoiSMF} second remark of page 278 see \S\ref{CorrectionMoi} below).
\begin{definition} 
$K^0_{\rm{flat}}(M)$ is the algebraic $K$-theory group of flat complex vector bundles
over $M$ modulo exact sequences. For any element $[E,\nabla_{\!E}]\in K^0_{\rm{flat}}(M)$,
($\nabla_{\!E}$ is a flat connection on $E$), its image in $\overset*K{}^0_{c_{\rm{tot}}}(M)$
is $b([E,\nabla_{\!E}])=(E,\nabla_{\!E},0)$.
\end{definition}
Of course $b$ takes its values in ${\mathcal M}K^0_{c_{\rm{tot}}}(M)$.
\begin{lemma}
$b$ defines a morphism from $K^0_{\rm{flat}}(M)$ to $\overset*K{}^0_{c_{\rm{tot}}}(M)$.
\end{lemma}
\begin{proof}
For an exact sequence of flat vector bundles as
\[0\longrightarrow (E',\nabla_{\!E'})\overset i\longrightarrow(E,\nabla_{\!E})
\overset p\longrightarrow
(E'',\nabla_{\!E''})\longrightarrow 0\]
if $s\colon E\to E'$ is such that $s\circ i={\rm{Id}}_{E'}$, then $\widetilde
c_{\rm{tot}}\big(\nabla_{\!E},(s\oplus p)^*(\nabla_{\!E'}\oplus\nabla_{\!E''})\big)=0$
and consequently $\widehat c_{\rm{tot}}\big(\nabla_{\!E},(s\oplus p)^*(\nabla_{\!E'}\oplus\nabla_{\!E''})\big)$ vanishes too. This is a general fact for exact sequences which respect (flat or nonflat) connections and any conjugation invariant polynomial (see lemma 6 in \cite{MoiPartie1} or proof of theorem 3.5 in \cite{MoiSMF} and combine it with the fact that
any conjugation invariant polynomial vanishes on any strictly upper triangular matrix).
\end{proof}
\subsection{Borel-type class:}\label{Borel}
In this context of multiplicative transgression, there exists some
equivalent to the Borel class defined in \cite{MoiPartie1} \S3.3.3.

Consider some smooth complex vector bundle $E$ endowed with some connection
$\nabla$ and some hermitian metric $h$. One defines the adjoint connection
$\nabla^{*}$ to $\nabla$ with respect to $h$ by the following formula
\begin{equation}{\tt x}.(h(s,s'))=h(\nabla_{\tt x}s,s')+h(s,\nabla_{\tt x}^*s')
\end{equation}
valid for any local sections $s$ and $s'$ of $E$ and any tangent vector
${\tt x}$ to $M$, ${\tt x}.f$ being the derivative of the function $f$
along the vector ${\tt x}$. One also defines the connection $\nabla^{u}
=\frac12(\nabla+\nabla^{*})$.
\begin{definition}\label{Mélanie}
The conjugate of $(E,\nabla,\alpha)\in\overset*K{}^0_{c_{\rm{tot}}}(M)$
is $(E,\nabla^*,\overline\alpha)$.
\end{definition}
\begin{lemma} The class of $(E,\nabla^*,\overline\alpha)$ in $\overset*K{}^0_{c_{\rm{tot}}}(M)$
does not depend of the hermitian metric used to evaluate $\nabla^*$.

$\overset*c_{\rm{tot}}$ is real (in the sense that $\overset*c_{\rm{tot}}(E,\nabla^*,\overline\alpha)
=\overline{\overset*c_{\rm{tot}}(E,\nabla,\alpha)}$).
\end{lemma}
\begin{proof}
If $h^0$ and $h^1$ are two different hermitian metrics on $E$, then choose any path $(h_t)_{t\in[0,1]}$ of hermitian metrics on $E$ with expected $h_0$ and $h_1$, and define for any $t$ the bundle automorphism $g_t$ of $E$ by
\begin{equation}\label{quoi}
h_t(s,s')=h_0(s,g_ts')
\end{equation}
for any $t$, the adjoint $\nabla^{*t}$ of $\nabla$ with respect
to $h_t$ is equal to $g_t^{-1}\nabla^{*0}g_t$ so that $\nabla^{*1}$ and $\nabla^{*0}$ are gauge conjugate by the gauge $g_1$ which is isotopic
to $g_0={\rm{Id}}_E$.
Thus, from \eqref{[}, $\widetilde c_{\rm{tot}}(\nabla^{*0},\nabla^{*1})$
and $\widehat c_{\rm{tot}}(\nabla^{*0},\nabla^{*1})$
are both exact.
The independence of the class of $(E,\nabla^*,\overline\alpha)$ on the hermitian metric follows
\eqref{reldeq} (and the facts that special addition of an exact form coincides with ordinary addition,
and exact forms vanish in $\Lambda$).

The curvatures of $\nabla$ and $\nabla^{*}$ are skew-adjoint, so that
their total Chern class are mutually complex conjugates. The reality of $\overset*c_{\rm{tot}}$ follows.
\end{proof}
This conjugation property extends to Chern-Simons forms (see \eqref{CS})
\begin{equation}\label{complexconjug}\widetilde c_{\rm{tot}}(\nabla_{\!0}^*,\nabla_{\!1}^*)=
\overline{\widetilde c_{\rm{tot}}(\nabla_{\!0},\nabla_{\!1})}
\qquad{\text{ and }}\qquad\widehat c_{\rm{tot}}(\nabla_{\!0}^*,\nabla_{\!1}^{*})=
\overline{\widehat c_{\rm{tot}}(\nabla_{\!0},\nabla_{\!1})}
\end{equation}
\begin{definition}\label{definitiondelaclassedeBorel}For any $(E,\nabla,\alpha)\in\overset*K{}^0_{c_{\rm{tot}}}(M)$ its Borel class is
\begin{equation}\label{defborel}
{\mathfrak B}^*(E,\nabla,\alpha)=
\widehat c_{\rm{tot}}(\nabla^*,\nabla)\stackrel\cap-\alpha\stackrel\cap+\overline\alpha\, \in\Lambda
\end{equation}
\end{definition}
\begin{lemma}${\mathfrak B}^*(E,\nabla,\alpha)$ is special imaginary (in the sense of \S\ref{complexe}).
It depends on the hermitian metric used to
calculate $\nabla^*$ only throw addition of an exact form.

${\mathfrak B}^*$ is a morphism from
$\overset*K{}^0_{c_{\rm{tot}}}(M)$ to $\Lambda$ (endowed with $\stackrel\cap+$).
\end{lemma}
\begin{proof}
Of course $\nabla^{u*}=\nabla^u$; thus $\widehat c_{\rm{tot}}(\nabla^*,\nabla)
=\widehat c_{\rm{tot}}(\nabla^*,\nabla^{u*})\stackrel\cap+\widehat c_{\rm{tot}}(\nabla^u,\nabla)$ (see \eqref{multcocycle}) is special imaginary because of \eqref{echange}.
${\mathfrak B}^*(E,\nabla,\alpha)$
is special imaginary
because special imaginary forms form a (special) additive subgroup for $\overset\cap+$.

If $f\colon E\overset\sim\longrightarrow F$ is a vector bundle isomorphism,
put $\nabla_{\!E}$ on $E$ and $\nabla_{\!F}$ on $F$ and some hermitian metric
$h^F$ on $F$. If $E$ is endowed with the pullback metric $f^*h^F$, then $f^*(\nabla_{\!F}^*)=(f^*\nabla_{\!F})^*$, so that from \eqref{multcocycle}
and \eqref{complexconjug}
\begin{align*}
\widehat c_{\rm{tot}}(\nabla_{\!E}^{*},\nabla_{\!E})&=
\widehat c_{\rm{tot}}\big(\nabla_{\!E}^{*},(f^*\nabla_{\!F})^*\big)\stackrel\cap+
\widehat c_{\rm{tot}}(\nabla_{\!F}^{*},\nabla_{\!F})\stackrel\cap+
\widehat c_{\rm{tot}}(f^*\nabla_{\!F},\nabla_{\!E})\\
{\mathfrak B}^*(E,\nabla_{\!E},\alpha)&=
{\mathfrak B}^*\big(F,\nabla_{\!F},\alpha\stackrel\cap+\widehat c_{\rm{tot}}
(\nabla_{\!E},f^*\nabla_{\!F})\big)
\end{align*}
Finally, the adjoint of a direct sum connection with respect to a direct sum hermitian metric
is the direct sum of the adjoint connections of each summand, the compatibility of ${\mathfrak B}^*$
with \eqref{relsum} then follows \eqref{multsum}.
\end{proof}
$\stackrel\cap-{\mathfrak B}^*$ is twice the imaginary part in $\overset*K{}^0_{c_{\rm{tot}}}$ since from \eqref{reldeq} and \eqref{defborel}:
\[(E,\nabla,\alpha)-\overline{(E,\nabla,\alpha)}=\iota(\stackrel\cap-{\mathfrak B}^*(E,\nabla,\alpha)\big)\, \quad{\text{ for any }}\ (E,\nabla,\alpha)\in\overset*K{}^0
_{c_{\rm{tot}}}(M)\]

The kernel of ${\mathfrak B}^*$ is the invariant subgroup of $\overset*K{}^0_{c_{\rm{tot}}}(M)$
under
the conjugation involution of definition \ref{Mélanie}. It is the subgroup generated by triples $(E,\nabla,\alpha)$ where $\nabla$ respects some hermitian
metric on $E$ and $\alpha$ is real.
\subsection{Nonfree multiplicative smooth $K$-theory:}\label{Malika}
Let ${\mathcal G}^-\subset\Omega^{\text{odd}}(M,{\mathbb C})$
be a subgroup for the operation $\stackrel\cap+$, and ${\mathcal G}^+\subset\Omega_{[+]}(M,{\mathbb C})$
be a subgroup (for the operation $\wedge$)
such that $1+d{\mathcal G}^-\subset{\mathcal G}^+$.
\begin{definition}
The (non free) multiplicative $K$-theory ${\mathcal M}K^0_{\mathcal G}(M)$ associated to ${\mathcal G}^\pm$ is the quotient by $\iota({\mathcal G}^-)$ of the subgroup $\overset*c{}_{\rm{tot}}^{-1}({\mathcal G}^+)$ of $\overset*K{}^0_{c_{\rm{tot}}}(M)$.
\end{definition}
The groups $\overset*K{}^0_{c_{\rm{tot}}}$ or
${\mathcal M}K^0_{c_{\rm{tot}}}$ (mentioned in \S\ref{jorba}) correspond to the particular
choices ${\mathcal G}^-=\{0\}$ and ${\mathcal G}^+=\Omega_{[+]}(M,{\mathbb C})$ or ${\mathcal G}^+=\{1\}$ respectively.

Let $\Lambda_{\mathcal G}$ be the quotient by ${\mathcal G}^-$ of the subgroup of $\Lambda$ (with respect to $\stackrel\cap+$)
consisting of forms $\alpha$ such that $1+d\alpha$ lies in ${\mathcal G}^+$.
\begin{proposition}\label{mel} The following sequence is exact
\[K^1_{\rm{top}}(M)\overset{Sc_{\rm{tot}}}\longrightarrow\Lambda_{\mathcal G}\overset\iota
\longrightarrow{\mathcal M}K^0_{\mathcal G}(M)\overset{\mathcal F}\longrightarrow K^0_{\rm{top}}(M)\overset{c_{\rm{tot}}}\longrightarrow
H^{\rm{even}}(M,{\mathbb C})\big/({\rm{Ker}}d\cap{\mathcal G}^+)\]
\end{proposition}
\begin{proof}
This sequence is clearly a complex, since elements of ${\mathcal M}K^0_{\mathcal G}(M)$
are constrained to have their total Chern class in ${\mathcal G}^+$.

Any element of $K^0_{\rm{top}}(M)$ with vanishing total Chern class in the quotient $H^{\rm{even}}(M,{\mathbb C})\big/({\rm{Ker}}d\cap{\mathcal G}^+)$ can be represented by formal difference of vector bundles
whose total Chern class lies in the image of ${\rm{Ker}}d\cap{\mathcal G}^+$
in $H^{\rm{even}}(M,{\mathbb C})$. It thus lies in ${\mathcal F}({\mathcal M}K^0_{\mathcal G}(M))$.

The exactness at $\Lambda_{\mathcal G}$ or ${\mathcal M}K^0_{\mathcal G}(M)$ is proved in the same way as
in the free case, using the fact that $\Lambda_{\mathcal G}$ and ${\mathcal M}K^0_{\mathcal G}(M)$
are quotiented by ${\mathcal G}^-$ and $\iota({\mathcal G}^-)$ respectively.
\end{proof}
The homotopy invariance of \cite{KaroubiTGCCS} \S4.8, expresses here in the following way:
\begin{lemma}\label{homo} Suppose that ${\mathcal G}^+\wedge{\mathcal G}^-\subset{\mathcal G}^-$. Set $\widetilde{\mathcal G}^+={\mathcal G}^+\oplus dt\wedge{\mathcal G}^-$ (this is a subset of
$\Omega(M\times[0,1],{\mathbb C})=\Omega(M,{\mathbb C})\oplus dt\wedge\Omega(M,{\mathbb C})$ which is a multiplicative subgroup of $\Omega_{[+]}(M,{\mathbb C})$). In the same way, choose some
additive subgroup ${\mathcal H}^+$ of $\Omega^{\rm{even}}(M,{\mathbb C})$
such that $d{\mathcal H}^+\subset{\mathcal G}^-$ and $d{\mathcal G}^-\wedge{\mathcal H}^+\subset{\mathcal H}^+$ (${\mathcal H}^+=\{0\}$ is always a possible choice) and set $\widetilde{\mathcal G}^-={\mathcal G}^-\oplus dt\wedge{\mathcal H}^+$, then
\[{\mathcal M}K^0_{\widetilde{\mathcal G}}(M\times[0,1])\cong{\mathcal M}K^0_{\mathcal G}(M)\]
\end{lemma}
Thus ${\mathcal M}K^0_{c_{\rm{tot}}}$ is homotopy invariant, while it is not the case of $\overset*K{}^0_{c_{\rm{tot}}}$.
\begin{proof}Let's first prove that $\Lambda_{\mathcal G}$
is isomorphic to $\Lambda_{\widetilde{\mathcal G}}$ (by pull-back of differential forms by the projection $M\times[0,1]\longrightarrow M$). The choices of $\widetilde{\mathcal G}^\pm$ ensure that this pullback is well defined and injective. Now if $\alpha_t+dt\wedge\beta_t$
is such that $\alpha_0=0$ and
\[1+d(\alpha_t+dt\wedge\beta_t)=1+d\alpha_t+dt\wedge\left(\frac{\partial\alpha_t}{\partial t}-d\beta_t\right)\ \in\ \widetilde{\mathcal G}^+\]
then $\gamma_t=\frac{\partial\alpha_t}{\partial t}-d\beta_t\, \in{\mathcal G}^-$ for any $t$, and one has
\[\alpha_t+dt\wedge\beta_t-d\left(\int_0^t\!\beta_sds\right)=\int_0^t\!
\alpha_sds+dt\wedge\beta_t-dt\wedge\beta_t-\int_0^t\!d\beta_sds=\int_0^t\!\gamma_sds\]
Thus $\alpha_t+dt\wedge\beta_t\, \in\widetilde{\mathcal G}^-+{\rm{Im}} d$ and the surjectivity of the pullback follows.

Similar arguments work for the case of $H^{\rm{even}}(M,{\mathbb C})\big/({\rm{Ker}}d\cap{\mathcal G}^+)$, and we can conclude using
proposition \ref{mel} and the five lemma.
\end{proof}
\subsection{Former version (from \cite{MoiSMF}):}\label{CorrectionMoi}
The version of this theory used in \cite{MoiSMF} (last remark of \S4 in page 278) corresponds to
the image of $\overset*K{}^0_{c_{\rm{tot}}}(M)$ by this
isomorphism:
\[(E,\nabla,\alpha)\in\overset*K{}^0_{c_{\rm{tot}}}(M)\, \longmapsto\, [E,\nabla,-c_{\rm{tot}}(\nabla^2)\wedge(\stackrel\cap-\alpha)]\]
The corresponding group with objects $[E,\nabla,\beta]$
(where $\beta\in\Lambda$ as $\alpha$), has the following relations, corresponding to \eqref{reldeq} and \eqref{relsum} respectively:
\begin{equation}\label{newrel}
\begin{aligned}{\text{if }}\quad f\colon E&\overset\sim\longrightarrow F
\quad{\text{ then }}\quad [E,\nabla_{\!E},\gamma+\widetilde c_{\rm{tot}}(f^*\nabla_{\!F},\nabla_{\!E})]=[F,\nabla_{\!F},\gamma]\\
[E,\nabla_{\!E},\beta]&+[G,\nabla_{\!G},\delta]=\\&=[E\oplus G,
\nabla_{\!E}\oplus\nabla_{\!G},\beta\wedge c_{\rm{tot}}(\nabla_{\!E}^2)+\delta\wedge c_{\rm{tot}}(\nabla_{\!G}^2)
+\beta\wedge d\delta]
\end{aligned}\end{equation}
(the second equality here should replace the erroneous formula just before the diagram in page 278 of \cite{MoiSMF}).
The map
$\iota$ is here replaced by
\begin{equation}\label{newiota}
\gamma\in\Lambda\longmapsto\big[E,\nabla,(1+d\gamma)^{-1}\big(\beta+\gamma\wedge c_{\rm{tot}}(\nabla^2)\big)\big]
-[E,\nabla,\beta]
\end{equation}
(this should also replace the erroneous definition given four lines
after the diagram in page 278 of \cite{MoiSMF}). And the differential form valued total Chern class reads now:
\begin{equation}\label{newChern}
c_{\rm{tot}}([E,\nabla,\beta])\, =\, c_{\rm{tot}}(\nabla^2)-d\beta
\end{equation}
The first line of \eqref{newrel} and \eqref{newChern} prove that
this version of the theory is eligible as a ``unstable''
(but ``free'') multiplicative $K$-theory in the sense of \cite{KaroubiTGCCS} \S4, the second line of \eqref{newrel}
and \eqref{newiota} (which is compatible to the relation used to Karoubi because it is there used with trivial $E$ with trivial connection $d$ and for differential forms $\gamma$ which are closed modulo the used subcomplex of the de Rham complex of $M$,
see \cite{KaroubiTGCCS} DEFINITION 4.2) make it to a ``stable''
(but ``free'')
multiplicative $K$-theory in the sense of \cite{KaroubiTGCCS}
\S5 with multiplicative total Chern class \eqref{newChern}.

The ``nonfree'' theory can then be defined in
exactly the same way as here before; in the case considered in \cite{MoiSMF}, ${\mathcal G}^+$ and ${\mathcal G}^-$ are both direct sums of forms of type $(p,q)$ such that $p\geq q$.

See \S\ref4 below for the total Chern class with values in Hodge-Deligne cohomology sketched in \cite{MoiSMF}.
\subsection{Multiplicative versus additive smooth $K$-theory:}
Here I recall the definition of the ``additive'' smooth $K^0$-theory (extended Bunke-Schick's smooth $K^0$-theory \cite{BunkeSchick}) as presented in \cite{MoiPartie1}: as in the multiplicative case, objects
are triples of the form $(E,\nabla_{\!E},\alpha)$ where $\nabla_{\!E}$ is a connection on the complex vector bundle $E$ on $M$, and $\alpha\in\Lambda$. If $f\colon E\overset\sim\longrightarrow F$ is a
vector bundle isomorphism, then $(E,\nabla_{\!E},\alpha)$ and $\big(F,\nabla_{\!F},\alpha+\widetilde{\rm{ch}}(\nabla_{\!E},\nabla_{\!F})\big)$ are equivalent.
\begin{definition*}$\widehat K^0_{\rm{ch}}(M)$ is the quotient of the free abelian group generated by such equivalence classes of triples, modulo the folowing relation :
\[(E,\nabla_{\!E},\alpha)+(G,\nabla_{\!G},\beta)=(E\oplus G,\nabla_{\!E}\oplus\nabla_{\!G},\alpha+\beta)\]
\end{definition*}
The counterpart of proposition \ref{mel} is the exactness of the following sequence:
\[K^1_{\rm{top}}(M)\overset{S{\rm{ch}}}\longrightarrow\Lambda\overset a
\longrightarrow \widehat K^0_{\rm{ch}}(M)\longrightarrow K^0_{\rm{top}}(M)\longrightarrow0\]
where the last morphism is the obvious forgetful map, $S{\rm{ch}}$ is the suspension of the Chern character on $K^1_{\rm{top}}$ constructed as in
the paragraph preceding \eqref{Patricia}, and $a$ is defined by $a(\alpha)=(E,\nabla,\alpha)-(E,\nabla,0)$ (for any $E$ and $\nabla$, $a$ is the notation of \cite{BunkeSchick}, the notation used in \cite{KaroubiTGCCS} \S5.3 or \cite{KaroubiCCFFHA} \S1.13 is $v$); of course $\Lambda$ is here endowed with its usual addition $+$ (and not $\stackrel\cap+$).

The differential form Chern character and the Borel class on $\widehat K^0_{\rm{ch}}(M)$
are given by ${\rm{ch}}(E,\nabla,\alpha)={\rm{ch}}(\nabla^2)-d\alpha$
and ${\mathfrak B}(E,\nabla,\alpha)=\widetilde{\rm{ch}}(\nabla^*,\nabla)
-\alpha+\overline\alpha$; the map $[E,\nabla]
\longmapsto(E,\nabla,0)$ defines a morphism $K^0_{\rm{flat}}(M)\longrightarrow\widehat K_{\rm{ch}}(M)$.

Remember the definition of $F$ from \eqref{Christine} and the definition of $\phi$ from just before theorem \ref{maitressePS}, and consider the maps
\begin{equation}\label{Dale}\begin{aligned}\psi_F\colon\alpha\in\Lambda&\longmapsto
\phi(\alpha)\wedge F\big(\phi(d\alpha)\big)\in\Lambda\\
\overset*\psi_F\colon(E,\nabla,\alpha)\in\widehat K^0_{\rm{ch}}(M)&\longmapsto(E,\nabla,\psi_F(\alpha)\big)\in\overset*K{}^0_{c_{\rm{tot}}}(M)
\end{aligned}\end{equation}
Lemma \ref{Christine} and the fact that $\phi$ is an automorphism of $\Lambda$ endowed with $+$ and commutes with $d$ when restricted to odd degree forms, prove that $\psi_F$ is an isomorphism from $\Lambda$ endowed with $+$ to $\Lambda$ endowed with $\stackrel\cap+$ whose inverse is
\[\psi_F^{-1}\colon\beta\longmapsto\phi^{-1}\big(\beta\wedge G(d\beta)\big)\]
($G$ was defined in \eqref{Christine}). The following result is
a consequence of theorem \ref{maitressePS}
\begin{theorem}\label{Virginia}$\overset*\psi_F$ is an isomorphism, and the following diagrams commute
\[\begin{CD}
K^1_{\rm{top}}(M)@>S{\rm{ch}}>>\Lambda @>a>>\widehat K^0_{\rm{ch}}(M)@>>>K^0_{\rm{top}}(M)\\
@V\Vert VV@V\psi_FVV@VV\overset*\psi_FV@VV\Vert V\\
K^1_{\rm{top}}(M)@>Sc_{\rm{tot}}>>\Lambda @>\iota>>\overset*K{}^0_{c_{\rm{tot}}}(M)@>>>K^0_{\rm{top}}(M)
\end{CD}\]
\[\begin{CD}\Omega^{\rm{even}}(M,{\mathbb C})@<{\rm{ch}}<<\widehat K^0_{\rm{ch}}(M)@>{\mathfrak B}>>\Lambda\\
@V\exp\circ\phi VV
@VV\overset*\psi_FV@VV\psi_FV\\\Omega_{[+]}(M,{\mathbb C})@<\overset*c_{\rm{tot}}<<\overset*K{}^0_{c_{\rm{tot}}}(M)
@>{\mathfrak B}^*>>\Lambda\end{CD}\qquad\qquad 
\begin{CD}K^0_{\rm{flat}}(M)@>>>\widehat K^0_{\rm{ch}}(M)\\
@V\Vert VV
@VV\overset*\psi_FV\\
K^0_{\rm{flat}}(M)@>b>>\overset*K{}^0_{c_{\rm{tot}}}(M)
\end{CD}
\]
\end{theorem}
(For the square with ${\mathfrak B}$, the point is that $\psi_F(\overline\alpha)=\overline{\psi_F(\alpha)}$. For the square with $\overset*c_{\rm{tot}}$, the point is that $1+d(\psi_F(\alpha))
=\exp\circ\phi(d\alpha)$).

The advantage of $\widehat K^0_{\rm{ch}}$ is that it can be endowed with an easily expressed product structure
(see \cite{BunkeSchick} definition 4.1)
\[(E,\nabla_{\!E},\alpha)\cup(E',\nabla_{\!E'},\alpha')=(E\otimes F,\nabla_{\!E}\otimes\nabla_{\!F},{\rm{ch}}(\nabla_{\!E}^2)\wedge\alpha'+\alpha\wedge{\rm{ch}}(\nabla_{\!F}^2)-\alpha\wedge d\alpha')\]
(this is much less easy for $\overset*K{}^0_{c_{\rm{tot}}}$!).
However, $\overset*K{}^0_{c_{\rm{tot}}}$
can be endowed with a total Chern class with values in integral Cheeger-Simons differential characters (see the next section). In fact,
the denominators in the Chern character limit its target space to rational differential characters only: see \cite{BunkeSchick} \S6.
This is the reason why Karoubi used the total Chern class in \cite{KaroubiCCFFHA} \S7.19 to obtain a class in nonrationalised Hodge-Deligne cohomology (see \S\ref4 below). The ``calcul symbolique'' he alludes to
corresponds to the above defined isomorphism $\overset*\psi_F$.
Another drawback of $\overset*K_{c_{\rm{tot}}}$
is that it only works (until now) at the $K^0$ level.

Finally, ``nonfree'' additive smooth $K^0$-theory $\, {\rm{ch}}^{-1}({\mathcal F}^+)\big/a({\mathcal F}^-)$ can be constructed from additive subgroups ${\mathcal F}^+$ of $\Omega^{\rm{even}}(M,{\mathbb C})$ and ${\mathcal F}^-$ of
$\Lambda$, such that ${\mathbb Z}\oplus d{\mathcal F}^-\subset{\mathcal F}^+$. 
In this case, $\overset*\psi_F$ induces an isomorphism
from this group to ${\mathcal M}K^0_{\mathcal G}$
with ${\mathcal G}^-=\psi_F({\mathcal F}^-)$
and ${\mathcal G}^+=\exp\circ\phi({\mathcal F}^+)$.
(The multiplicative structure of $\widehat K_{\rm{ch}}(M)$
descends to $\, {\rm{ch}}^{-1}({\mathcal F}^+)\big/a({\mathcal F}^-)$ under the extra hypotheses
${\mathcal F}^+\wedge{\mathcal F}^+\subset{\mathcal F}^+$ and ${\mathcal F}^-\wedge{\mathcal F}^+\subset{\mathcal F}^-$).
\section{Differential character valued Chern class:}
\subsection{Cheeger-Simons differential characters:}
Cheeger and Simons' differential characters $\widehat H^\bullet(M,{\mathbb R}/{\mathbb Z})$ were defined in \cite{CheegerSimons}
as the singular cocycles with values in ${\mathbb R}/{\mathbb Z}$
whose restrictions to boundaries are ${\mathbb R}/{\mathbb Z}$ reductions
of cocycles lying in the image of real valued
$C^\infty$ differential forms by integration on chains.

We will consider ${\mathbb C}/{\mathbb Z}$ valued objects here.
Let $\Lambda^{\rm{even}}_0$ denote the space of closed even-degree differential forms whose de Rham cohomology class is integral,
the two first exact sequences of \cite{CheegerSimons} theorem 1.1 read here
\begin{equation}\label{bleble}
\begin{aligned}
0\longrightarrow H^{\rm{odd}}({\mathbb C}/{\mathbb Z})
\overset{\iota_1}\longrightarrow &\widehat H^{\rm{odd}}(M,{\mathbb C}/{\mathbb Z})
\overset{\delta_1}\longrightarrow\Lambda^{\rm{even}}_0\longrightarrow0\\
0\longrightarrow \Lambda/H^{\rm{odd}}(M,{\mathbb Z})
\overset{\iota_2}\longrightarrow &\widehat H^{\rm{odd}}(M,{\mathbb C}/{\mathbb Z})\overset{\delta_2}
\longrightarrow H^{\rm{even}}(M,{\mathbb Z})\longrightarrow0
\end{aligned}\end{equation}
Here the graduation convention on $\widehat H$ is the one of \cite{CheegerSimons}, in spite of the fact that the nowadays usual one is shifted by one with respect to it. In fact both lines are graduated, in the sense that they decompose to direct sum of lines with the same grading on the two left terms, and the same plus 1 on the right term.
In particular $\widehat H^{-1}(M,{\mathbb C}/{\mathbb Z})\cong H^0(M,{\mathbb Z})$. Formula 2) in corollary 1.2 of \cite{CheegerSimons} implies that $\delta_1\circ\iota_2=d$ is the usual exterior differential.

There is a multiplicative structure $*$ on differential characters,
whose unity is $1\in\widehat H^{-1}(M,{\mathbb C}/{\mathbb Z})
\cong H^0(M,{\mathbb Z})$ and for which $\delta_1$ and $\delta_2$ are ring homomorphisms. It is a consequence of formula (1.15) in \cite{CheegerSimons} that for any $\alpha\in\Lambda$ and any $f\in\widehat H^{\rm{odd}}(M,{\mathbb C}/{\mathbb Z})$, one has
\begin{equation}\label{IranzoSarah}\iota_2(\alpha)*f=\iota_2\big(\alpha\wedge\delta_1(f)\big)\end{equation}
In particular, $\Lambda/H^{\rm{odd}}(M,{\mathbb Z})$ is an ideal in $\widehat H^{\rm{odd}}(M,{\mathbb C}/{\mathbb Z})$ and for any $\alpha$ and $\beta\in\Lambda$:
\[\iota_2(\alpha)*\iota_2(\beta)=\iota_2(\alpha\wedge d\beta)=\iota_2(\beta\wedge d\alpha)\]
As a consequence of this, one obtains
\begin{equation}\label{plusiota}
\big(1+\iota_2(\alpha)\big)*\big(1+\iota_2(\beta)\big)=1+\iota_2(\alpha\stackrel\cap+\beta)
\end{equation}
In particular $\big(1+\iota_2(\alpha)\big)^{-1}=1+\iota_2\big(\stackrel\cap-\alpha\big)$.

The initial theory of differential characters \cite{CheegerSimons}
is a real theory $\widehat H^{\rm{odd}}(M,{\mathbb R}/{\mathbb Z})$
from which the complex theory is obtained by adding (as direct sum)
the space $i\, {\mathfrak{Im}}\Lambda$ of odd degree purely imaginary differential forms modulo exact forms:
\begin{equation}\label{PreBorel}
\widehat H^{\rm{odd}}(M,{\mathbb C}/{\mathbb Z})
=\widehat H^{\rm{odd}}(M,{\mathbb R}/{\mathbb Z})
\oplus i\, {\mathfrak{Im}}\Lambda
\end{equation}
$i\, {\mathfrak{Im}}\Lambda$ enters as a direct summand in the left term
of the second line of \eqref{bleble}, and as $H^{\rm{odd}}(M,{\mathbb C}/{\mathbb Z})\cong H^{\rm{odd}}(M,{\mathbb R}/{\mathbb Z})\oplus iH^{\rm{odd}}(M,{\mathbb R})$, it gives rise to the following
exact subsequence of the first line of \eqref{bleble}
\begin{equation}\label{imagexact}
0\longrightarrow iH^{\rm{odd}}(M,{\mathbb R})\longrightarrow
i\, {\mathfrak{Im}}\Lambda\overset d\longrightarrow i\, 
d\Omega^{\rm{odd}}(M,{\mathbb R})\longrightarrow0
\end{equation}
The decomposition \eqref{PreBorel} is considered as a decomposition in real part $\oplus$ imaginary part, and it gives rise to some obvious ``complex conjugation''
involution of $\widehat H^{\rm{odd}}(M,{\mathbb C}/{\mathbb Z})$ which will be classically denoted by $f\longmapsto\overline f$.

It follows \eqref{IranzoSarah} and a straightforward calculation
that the conjugation is compatible with the product $*$, and thus with the inverse. It is also obviously compatible with $\iota_1$, $\iota_2$ and $\delta_1$. Finally $\delta_2(f)=\delta_2(\overline f)$ for all $f$.
\subsection{Product decomposition:}
Consider now the subset $\widehat H_{[+]}(M,{\mathbb C}/{\mathbb Z})=\delta^{-1}_1\big(\Omega_{[+]}(M,{\mathbb C})\big)$ of
$\widehat H^{\rm{odd}}(M,{\mathbb C}/{\mathbb Z})$ consisting of elements of the form $1+f$ with $f$ of positive multidegree. It is a commutative group for the law $*$, and $\delta_1$ is a group morphism from it to $\Omega_{[+]}(M,{\mathbb C})$ (with group law $\wedge$). In the same way $\delta_1$ is a morphism from it to
$1+\mathop\oplus\limits_{k\geq1}H^{2k}(M,{\mathbb Z})$ with the cap product as group law.

The counterpart of the exact sequences \eqref{bleble} read
\begin{equation}\label{blebleble}
\begin{aligned}
0\longrightarrow H^{\rm{odd}}({\mathbb C}/{\mathbb Z})
\overset{1+\iota_1}{-\!\!\!\longrightarrow} &\widehat H_{[+]}(M,{\mathbb C}/{\mathbb Z})
\overset{\delta_1}\longrightarrow\Omega_{[+]}(M,{\mathbb C})\cap\Lambda_0^{\rm{even}}\longrightarrow\{1\}\\
0\longrightarrow \Lambda/H^{\rm{odd}}(M,{\mathbb Z})
\overset{1+\iota_2}{-\!\!\!\longrightarrow} &\widehat H_{[+]}(M,{\mathbb C}/{\mathbb Z})\overset{\delta_2}
\longrightarrow 1+\mathop\oplus\limits_{k\geq1}H^{2k}(M,{\mathbb Z})\longrightarrow\{1\}
\end{aligned}\end{equation}
where now $\Lambda/H^{\rm{odd}}(M,{\mathbb Z})$ is endowed with the group law $\overset\cap+$ (see \eqref{plusiota}).Denote by $\widehat H_{[+]}(M,{\mathbb R}/{\mathbb Z})=\delta^{-1}_1\big(\Omega_{[+]}(M,{\mathbb R})\big)$ the real elements in $\widehat H_{[+]}(M,{\mathbb C}/{\mathbb Z})$,
\begin{lemma}\label{decompo}Any element $f$ of $\widehat H_{[+]}(M,{\mathbb C}/{\mathbb Z})$
decomposes in a unique way as a product $\, f=g*\big(1+\iota_2(\beta)\big)$
where
$\beta\in\Lambda/H^{\rm{odd}}(M,{\mathbb Z})$ is special imaginary and $g\in\widehat H_{[+]}(M,{\mathbb R}/{\mathbb Z})$ is real.
\end{lemma}
In particular, it follows from \eqref{plusiota} that $1+\iota_2(\beta)$ is of modulus one in the sense that
$\big(1+\iota_2(\beta)\big)*\overline{\big(1+\iota_2(\beta)\big)}=1$ (using the complex conjugation from \eqref{PreBorel}
on $\widehat H^{\rm{odd}}(M,{\mathbb C}/{\mathbb Z})$), so that
the decomposition of the lemma corresponds to some decomposition of the same kind
as at the beginning of \S\ref{complexe}. However,
there may be some elements in $\widehat H_{[+]}(M,{\mathbb C}/{\mathbb Z})$ which
are simultaneously real and of modulus 1 (because of torsion in $H^\bullet(M,{\mathbb Z})$ and $H^\bullet(M,{\mathbb C}/{\mathbb Z})$) so
that a decomposition as a product of a real and a modulus $1$ element may be nonunique without the constraint that the modulus $1$ element be of the form $1+\iota_2(\beta)$.
\begin{proof}
Let $f\in\widehat H_{[+]}(M,{\mathbb C}/{\mathbb Z})$, then $\delta_1(f)$
decomposes uniquely as $\rho\wedge\theta$ with $\rho$ real and $\theta$ of modulus 1 in $\Omega_{[+]}(M,{\mathbb C})$. From the last statement of lemma \ref{SophieMusicale} and lemma \ref{sophie}, one obtains the existence of some special imaginary $\alpha\in\Lambda$
such that $1+d\alpha=\delta_1\big(1+\iota_2(\alpha)\big)=\theta$.
It then follows that $\rho=\delta_1\Big(f*\big(1+\iota_2(\alpha)\big)^{-1}\Big)$.

Denote by $\gamma\in\, i\, {\mathfrak{Im}}\Lambda$ the imaginary part
in the sense of decomposition \eqref{PreBorel} of $f*\big(1+\iota_2(\alpha)\big)^{-1}$. Of course $\gamma$ is a closed form (see \eqref{imagexact}) and so is $\gamma\wedge\rho^{-1}$. One then defines $\beta$ by $\beta=\alpha+\gamma\wedge\rho^{-1}=\alpha\stackrel\cap+\gamma\wedge\rho^{-1}$.
Any closed purely imaginary form is special imaginary in the sense of \S\ref{complexe}, so that $\beta$ is special imaginary itself.
Then using \eqref{IranzoSarah} and \eqref{plusiota}, one obtains the reality of
\begin{align*}
g&=f*\big(1+\iota_2(\beta)\big)^{-1}=f*\big(1+\iota_2(\alpha)\big)^{-1}
*\big(1+\iota_2(\gamma\wedge\rho^{-1})\big)^{-1}\\
&=f*\big(1+\iota_2(\alpha)\big)^{-1}
*\big(1+\iota_2(-\gamma\wedge\rho^{-1})\big)=f*\big(1+\iota_2(\alpha)\big)^{-1}+\iota_2(-\gamma)
\end{align*}
This proves the existence of the decomposition.

The uniqueness is due to the fact that any nonvanishing special imaginary form has a nonvanishing (nonspecial) imaginary part, and from the exact sequence \eqref{imagexact}, which together prove that no element of $\widehat H^{\rm{odd}}(M,{\mathbb C}/{\mathbb Z})$ can be simultaneously real, different from $1$, and of the form $1+\iota_2(\beta)$ with $\beta$ special imaginary.
\end{proof}
\subsection{Definition of the class:}
\begin{theorem}\label{Angelique}
There exists a unique total Chern class morphism
\[\overset\cup c_{\rm{tot}}\colon\overset*K{}^0_{c_{\rm{tot}}}\longrightarrow\widehat H_{[+]}(M,{\mathbb C}/{\mathbb Z})\]
with following compatibility properties:
\begin{equation}\label{jambe}\begin{aligned}\delta_1\big(\overset\cup c_{\rm{tot}}(E,\nabla,\alpha)\big)
&=\overset*c_{\rm{tot}}(E,\nabla,\alpha)\\
\delta_2\big(\overset\cup c_{\rm{tot}}(E,\nabla,\alpha)\big)
&=c_{\rm{tot}}([E])\\
\overset\cup c_{\rm{tot}}\big(\iota(\alpha)\big)&=\big(1+\iota_2(\alpha)\big)^{-1}
\end{aligned}\end{equation}
$\overset\cup c_{\rm{tot}}$ is real in the sense that it is compatible with the ``complex'' conjugations defined in definition \ref{Mélanie} and just after \eqref{imagexact}. Its relation with ${\mathfrak B}^*$ is given by the decomposition of $\overset\cup c_{\rm{tot}}(E,\nabla,\alpha)$ as in lemma \ref{decompo}:
\[\overset\cup c_{\rm{tot}}(E,\nabla,\alpha)=
f*\left(1+\iota_2\big({}^{1}_{\overline2}{\mathfrak B}^*(E,\nabla,\alpha)\big)\right)\qquad{\text{where }}\quad f\in\widehat H^{\rm{odd}}(M,{\mathbb R}/{\mathbb Z})\]
\end{theorem}
This theorem answers the question raised by U. Bunke in \cite{BunkeQuestionOberwolfach}.
\begin{proof}
There is a theory of characteristic classes for vector bundles with connections
with values in differential characters developped in \S2 (especially theorem 2.2) and \S4 (for the total 
Chern class, which is called there ``Chern character'') of \cite{CheegerSimons}.
However, there is some ambiguity in \cite{CheegerSimons} because the authors use
real differential characters, and seem to construct a class for any connection.
In particular, they do not precise which classifying space (\cite{NaRam1} or \cite{NaRam2}) they use. It is possible to use \cite{NaRam1} only for
connections which respect some hermitian metric, and in this case, the Chern-Weil total Chern class is a real form, and the obtained differential character total Chern class lies in $\widehat H^{\rm{odd}}(M,{\mathbb R}/{\mathbb Z})$.
For general connections, it is necessary to use \cite{NaRam2}, the Chern-Weil
total chern class need not be real, and it is thus impossible to obtain a real differential character valued total Chern class. However, the technique explained in
\cite{CheegerSimons} Theorem 2.2 works, and provides a total Chern class
with values in $\widehat H^{\rm{odd}}(M,{\mathbb C}/{\mathbb Z})$ (in $\widehat H_{[+]}(M,{\mathbb C}/{\mathbb Z})$ in fact), which will be denoted by $\overset\vee c_{\rm{tot}}$ here.
Of course the two obtained classes agree for connections which respect some hermitian metrics because of universality.
The two first statements of \cite{CheegerSimons} theorem 2.2 read here
\begin{equation}\label{bras}\delta_1\big(\overset\vee c_{\rm{tot}}(E,\nabla)\big)
=c_{\rm{tot}}(\nabla^2)\qquad{\text{and}}\qquad
\delta_2\big(\overset\vee c_{\rm{tot}}(E,\nabla)\big)
=c_{\rm{tot}}([E])\end{equation}

Moreover, it follows from \cite{CheegerSimons}
proposition 2.9 that if $\nabla_{\!0}$ and $\nabla_{\!1}$ are two connections on the same vector bundle $E$, then
\begin{equation}\label{anomalyvee}
\overset\vee c_{\rm{tot}}(E,\nabla_{\!1})=\overset\vee c_{\rm{tot}}(E,\nabla_{\!0})+\iota_2\big(\widetilde c_{\rm{tot}}(\nabla_{\!0},\nabla_{\!1})\big)
\end{equation}
Combining this with \eqref{IranzoSarah} yields
\[\overset\vee c_{\rm{tot}}(E,\nabla_{\!1})*\big(\overset\vee c_{\rm{tot}}(E,\nabla_{\!0})\big)^{-1}=1+\iota_2\big(\widehat c_{\rm{tot}}(\nabla_{\!0},\nabla_{\!1})\big)\]

From this, \eqref{reldeq} and \eqref{plusiota}, one obtains a map
\begin{equation}\label{defcup}
(E,\nabla,\alpha)\in\overset*K{}^0_{c_{\rm{tot}}}(M)\, \longmapsto\, \overset\vee c_{\rm{tot}}(E,\nabla)*\big(1+\iota_2(\alpha)\big)^{-1}\in\widehat H_{[+]}(M,{\mathbb C}/{\mathbb Z})
\end{equation}
which will be denoted by $\overset\cup c_{\rm{tot}}$. \eqref{relsum}, \eqref{plusiota} and the multiplicativity of $\overset\vee c_{\rm{tot}}$
(proved in \cite{CheegerSimons} theorems 4.6 and 4.7) prove that $\overset\cup c_{\rm{tot}}$ is a group morphism. The compatibility properties \eqref{jambe} follow definition \ref{italiénisante}, \eqref{bras}, \eqref{defcup} and the facts that $\delta_1\circ\iota_2=d$ and that $\delta_1$ is a ring homomorphism.

The unicity of $\overset\cup c_{\rm{tot}}$ can be proved in exactly the same way as the unicity of the rationalised differential character valued Chern character
on $\widehat K^0_{\rm{ch}}(M)$
proved by U. Bunke and Th. Schick in \cite{BunkeSchick} \S\S6.2.2 and 6.2.3.

Take now some element $(E,\nabla,\alpha)\in\overset*K{}^0_{c_{\rm{tot}}}(M)$.
Put on $E$ any hermitian metric $h$, and consider the adjoint transpose $\nabla^*$ and the unitary connection $\nabla^u$
obtained from $\nabla$ and $h$ as just before definition \ref{Mélanie}.
Consider the unique real forms $\beta$ and $\delta$ and special imaginary forms
$\gamma$ and $\zeta$ entering in the decompositions
of $\alpha$ and $\widehat
c_{\rm{tot}}(\nabla^u,\nabla)$ as in lemma \ref{deCadix}:
\[\alpha=\beta\stackrel\cap+\gamma\qquad{\text{ and }}\qquad\widehat
c_{\rm{tot}}(\nabla^u,\nabla)=\delta\stackrel\cap+\zeta\]
It follows from \eqref{multcocycle}, \eqref{echange}, \eqref{complexconjug} and \eqref{defborel} that:
\begin{equation}\label{rraallbbooll}\widehat
c_{\rm{tot}}(\nabla^u,\nabla^*)=\delta\stackrel\cap-\zeta\qquad{\text{ and }}\qquad{\mathfrak B}^*(E,\nabla,\alpha)=\zeta\stackrel\cap+\zeta\stackrel\cap-\gamma\stackrel\cap-\gamma\end{equation}

Then, using \eqref{plusiota}, \eqref{anomalyvee} and \eqref{defcup}, one obtains
\begin{align*}
\overset\cup c_{\rm{tot}}(E,\nabla,\alpha)&=\overset\cup c_{\rm{tot}}(E,\nabla^u,0)
*\big(1+\iota_2(\beta\stackrel\cap+\gamma)\big)^{-1}*\big(1+\iota_2(\widehat c_{\rm{tot}}(\nabla^u,\nabla))\big)\\
&=\overset\vee c_{\rm{tot}}(E,\nabla^u)*\big(1+\iota_2(\delta\stackrel\cap-\beta)\big)*\big(1+\iota_2(\zeta\stackrel\cap-\gamma)\big)
\end{align*}
The last statement of the theorem follows this last equation, \eqref{rraallbbooll} and the fact that $\overset\vee c_{\rm{tot}}(E,\nabla^u)$ and $\delta\stackrel\cap-\beta$ are both real. The reality of $\overset\cap c_{\rm{tot}}$ is a consequence of the facts that the same calculation for $\overset\cup c_{\rm{tot}}(E,\nabla^*,\overline\alpha)$ provides
\[\overset\cup c_{\rm{tot}}(E,\nabla^*,\overline\alpha)=\overset\vee c_{\rm{tot}}(E,\nabla^u)*\big(1+\iota_2(\delta\stackrel\cap-\beta)\big)*\big(1+\iota_2(\stackrel\cap-\zeta\stackrel\cap+\gamma)\big)\]
and that the inverse of the modulus one element $1+\iota_2(\zeta\stackrel\cap-\gamma)$ coincides with its conjugate $1+\iota_2(\stackrel\cap-\zeta\stackrel\cap+\gamma)$.
\end{proof}
\subsection{Some extra properties of $\overset\cup c_{\rm{tot}}$:}\label4
For any integer $n$, let ${\mathbb C}^n$ be the trivial rank $n$ vector bundle, and $d$ its canonical connection. The map $n\longmapsto ({\mathbb C}^n,d,0)$ extends to a morphism ${\mathbb Z}\longrightarrow\overset*K{}^0_{c_{\rm{tot}}}(M)$.
\begin{lemma*}
$\overset\cup c_{\rm{tot}}$ induces an isomorphism \[\big(\overset*K{}^0_{c_{\rm{tot}}}(M)\big/{\mathbb Z}\big)\otimes{\mathbb Q}\overset\sim\longrightarrow\widehat H_{[+]}(M,{\mathbb C}/{\mathbb Z})\otimes{\mathbb Q}\]
\end{lemma*}
\begin{proof}
It is a direct consequence of the five lemma, the exact sequences \eqref{pasinvarhomot} and \eqref{blebleble} (the second one), the compatibility relations \eqref{jambe} and the facts that $c_{\rm{tot}}$ induces an isomorphism between $\left(K^0_{\rm{top}}(M)\big/{\mathbb Z}\right)\otimes{\mathbb Q}$ and $1+\mathop\oplus\limits_{k\geq1}H^{2k}(M,{\mathbb Q})$ and $Sc_{\rm{tot}}$ induces an isomorphism between $K^1_{\rm{top}}(M)\otimes{\mathbb Q}$ and $H^{odd}(M,{\mathbb Q})$.
\end{proof}
There is a morphism $\ln\colon\widehat H_{[+]}(M,{\mathbb C}/{\mathbb Z})\longrightarrow\widehat H^{\rm{odd}}(M,{\mathbb C}/{\mathbb Q})$ (there are denominators in the series $\ln(1+X)$). It follows \eqref{IranzoSarah} and the fact that $\delta_1\circ\iota_2=d$ that for any $\alpha\in\Lambda$:
\begin{equation}\label{Armel}\ln(1+\iota_2(\alpha))=\iota_2(\alpha\wedge G(d\alpha))\end{equation}
Let $\phi$ be defined on $\widehat H^{\rm{odd}}(M,{\mathbb C}/{\mathbb Z})$
as it is on differential forms (see just before theorem \ref{maitressePS}), namely
$\phi$ multiplies the $(2k-1)$-degree part by $(-1)^{k-1}(k-1)!$ for all $k$
(and vanishes on the $-1$ degree part, corresponding to ${\mathbb Z}$).
The values taken by the logarithm in $\widehat H^{\rm{odd}}(M,{\mathbb C}/{\mathbb Z})$ have all vanishing $-1$ degree parts, so that one can apply $\phi^{-1}$ to them. \eqref{Armel} and the formula $c_{\rm{tot}}=\exp\circ\phi({\rm{ch}})$ yield
\begin{theorem}\label{Karen} The following diagram commutes:
\[\begin{CD}\overset*K{}^0_{c_{\rm{tot}}}(M)@>\overset*\psi{}_{F}^{-1}>>\widehat K^0_{\rm{ch}}(M)\\@V\overset\cup c_{\rm{tot}}VV@VV\widehat{\rm{ch}}-{\rm{rk}}V\\
\widehat H_{[+]}(M,{\mathbb C}/{\mathbb Z})@>\phi^{-1}\circ\ln>>\widehat H^{\rm{odd}}(M,{\mathbb C}/{\mathbb Q})
\end{CD}\]
where $\widehat{\rm{ch}}-{\rm{rk}}$ is the difference of the Chern character
of \cite{BunkeSchick} \S6 and the rank.
\end{theorem}

Finally, taking ${\mathcal G}^+$ and ${\mathcal G}^-$ as in \S\ref{Malika},
one can define a ``multiplicative cohomology'' group as ${\mathcal M}H_{\mathcal G}(M)=\delta_1^{-1}({\mathcal G}^++{\rm{Im}}d)\big/(1+\iota_2({\mathcal G}^-))$, which enters in an exact sequence of the same type as in proposition \ref{mel}, and turns out to be a smooth homotopy invariant in the same conditions (and the same hypotheses on ${\mathcal G}^\pm$) as ${\mathcal M}K^0_{\mathcal G}(M)$ (see lemma \ref{homo}).
The total Chern class $\overset\cup c_{\rm{tot}}$ then provides a morphism ${\mathcal M}K^0_{\mathcal G}(M)\longrightarrow{\mathcal M}H_{\mathcal G}(M)$. If $M$ is a compact K\"ahler manifold, take
${\mathcal G}^+$ and ${\mathcal G}^-$ to be the sums of (respectively even and odd degree) differential forms of type $(p,q)$ with $p\geq q$ and with degree $0$ component equal to $1$ for ${\mathcal G}^+$.

Bearing in mind that here the total Chern class only on degree $0$ $K$-theory is treated, this provides
the ``multiplicative'' counterpart of the example given in \cite{KaroubiCCFFHA} \S\S7.12-19 (and also \cite{Brylinski} lemma 2, where it appears as a quotient of a subgroup of ``restricted'' differential characters, or \cite{Felisatti} formula just before (8) and theorem 2.3; there are corresponding constructions in \cite{Zucker}. Some of these constructions generalise to the quasiprojective or noncompact analytic cases, which are not considered here because I do not want to consider logarithmic forms on noncompact manifolds).
The obtained ${\mathcal M}H_{\mathcal G}$ coincides with the multiplicative subgroup of Hodge-Deligne cohomology consisting of objects whose $0$-degree component equals $1$, and the morphism induced by $\overset\cup c_{\rm{tot}}$
(composed with the morphism $\overset*\psi_F$ of \S\ref{Dale}) coincides with the total Chern class considered in \cite{KaroubiCCFFHA} \S\S7.18 and 7.19.


\begin{thebibliography}9
\bibitem{MoiPartie1} A. Berthomieu: {\emph{Direct image for relative and multiplicative $K$-theories from transgression of the families index theorem, part 1.}},
preprint at {\tt{arXiv:math.DG/0611281}}
\bibitem{MoiPartie2et3} A. Berthomieu: {\emph{Direct image for relative and multiplicative $K$-theories from transgression of the families index theorem, part 2. and 3.}},
preprints at {\tt{arXiv:math.DG/0703916}} and {\tt{arXiv:0804/0728}}.
\bibitem{MoiSMF} A. Berthomieu: {\emph{Proof of Nadel's conjecture and direct image for relative $K$-theory}}, Bull. Soc. Math. France {\bf{130}} (2), 2002, pp. 253-307. 
\bibitem{Brylinski} J.-L. Brylinski: {\emph{Comparison of the Beilinson-Chern classes with the Chern-Cheeger-Simons classes}} in {\emph{Advances in Geometry}}, Progr. Math. {\bf{172}}, Birkh\"auser Boston, Boston, MA (1999) 95-105.
\bibitem{BunkeQuestionOberwolfach} U. Bunke: question in the problem session, Oberwolfach Reports, Vol. {\bf{3}}, Nr 1 (2006) page 797.
\bibitem{BunkeSchick} U. Bunke, Th. Schick: {\emph{Smooth $K$-theory}} preprint at {\tt{arXiv:0707/0046}}
\bibitem{CheegerSimons} J. Cheeger, J. Simons: {\emph{Differential characters and geometric invariants}}, Lecture Notes in Math. 1167 Springer, New York, 1985, pp. 50-80.
\bibitem{Felisatti} M. Felisatti: {\emph{Differential characters and mulyiplicative cohomology}}, $K$-theory {\bf{18}} (1999) 267-276.
\bibitem{Hopkins} M. J. Hopkins, I. M. Singer: {\emph{Quadratic functions in geometry, topology and $M$-theory}}, J. of Diff. Geom. {\bf{70}} n$^\circ$3
(2005) 329-452.
\bibitem{KaroubiTGCCS} M. Karoubi: {\emph{Th\'eorie g\'en\'erale des classes caract\'eristiques secondaires}}, $K$-theory {\bf{4}}, 1990, pp. 55-87.
\bibitem{KaroubiCCFFHA} M. Karoubi: {\emph{Classes caract\'eristiques de fibr\'es feuillet\'es, holomorphes ou alg\'ebriques}}, $K$-theory {\bf{8}}, 1994, pp. 153-211.
\bibitem{NaRam1} M. S. Narasimhan, S. Ramanan: {\emph{Existence of universal connections}}, Amer. J. of Math. {\bf{83}} (1961) 563-572.
\bibitem{NaRam2} M. S. Narasimhan, S. Ramanan: {\emph{Existence of universal connections II}}, Amer. J. of Math. {\bf{85}} (1963) 223-231. 
\bibitem{Zucker}S. Zucker: {\emph{The Cheeger-Simons invariant as a Chern class}}, in {emph{Algebraic analysis, geometry and number theory (Baltimore, MD, 1988)}}, John Hopkins Univ. Press, Baltimore, MD (1989) 397-417.
\end{thebibliography}
\end{document}